\documentclass{amsart}
\usepackage{graphicx}
\usepackage{amsmath}
\usepackage{color}

\newtheorem{thm}[equation]{Theorem}
\newtheorem{pro}[equation]{Proposition}
\newtheorem{cor}[equation]{Corollary}
\newtheorem{lem}[equation]{Lemma}

\theoremstyle{definition}
\newtheorem{exa}[equation]{Example}

\newtheorem{DEF}[equation]{Definition}
\newtheorem{rem}[equation]{Remark}

\usepackage{amsfonts,amssymb,amscd,amsmath,enumerate,verbatim,calc}

\def\a{\alpha}
\def\aa{\mathcal A}

\def\b{\beta}

\def\e{\epsilon}

\def\bbbf{\mathbb{F}}

\def\LL{\mathcal{L}}

\def\lam{\lambda}
\def\Lam{\Lambda}

\def\hh{\mathcal H}

\def\v{{\mathcal V}}

\def\bbbz{{\mathbb Z}}

\def\bbbk{{\mathbb K}}

\def\lm{\lambda}
\def\e{E}

\def\g{G}

\def\vv{{\mathcal V}}
\def\uu{{\mathcal U}}
\def\aa{{\mathcal A}}

\def\lam{\Lambda_I^{\bar{i}}}

\def\e{\epsilon}
\def\b{\beta}
\def\a{\alpha}
\def\lam{\lambda}
\def\Lam{\Lambda}
\begin{document}

\markboth{Split  $3$-Lie-Rinehart color algebras} {V. Khalili}

\date{}

\centerline{\bf Split   $3$-Lie-Rinehart color algebras }

\vspace{.5cm}\centerline{Valiollah Khalili\footnote[1]{Department
of mathematics, Faculty of sciences, Arak University, Arak 385156-8-8349, Po.Box: 879, Iran. 
 V-Khalili@araku.ac.ir\\
\hphantom{ppp}2000 Mathematics Subject Classification(s):17B05, 17B22, 17B60, 17A60.
\\
\hphantom{ppp}Keywords: $3-$Lie-Rinehart color algebra, split algebra, root space, weight space and structure theory. }\;\;}

\vspace{1cm} \noindent ABSTRACT:

In this paper we introduce a class of $3-$color algebras which are called split  $3-$Lie-Rinehart color algebras as the natural generalization of the one of split Lie-Rinehart algebras.  We characterize their inner structures by developing techniques of connections
 of root systems and weight systems associated to a splitting Cartan subalgebra. We show that such a tight  split  $3-$Lie-Rinehart color algebras $(\LL, A)$  decompose as the orthogonal
 direct sums  $\LL =\oplus_{i\in I}\LL_i$ and $A =\oplus_{j\in J}A_j,$
 where any $\LL_i$ is a non-zero graded ideal of $\LL$ stisfying $[\LL{i_1}, \LL{i_2}, \LL{i_3}]=0$ if $i_1, i_2, i_3\in I$ be different from each other and 
 any $A_j$ is a non-zero graded ideal of $A$  stisfying $A_{j_1}A_{j_2}=0$ if $J_1\neq j_2.$ Both decompositions satisfy that for any $i\in I$  there
 exists a unique $j\in J$ such that $A_j\LL_i = 0$. Furthermore, any $(\LL_i, A_j)$ is a split $3-$Lie-Rinehart color algebra. Also,  under certain conditions, it is shown that the above decompositions of $\LL$ and $A$ are by means of the family of their, respective, simple ideals.

\vspace{1cm} \setcounter{section}{0}
\section{Introduction}\label{introduction}
The notion of
$n-$Lie algebra  introduced by Filippov in 1985 (see \cite{F}), is a natural generalization of
the concept of a Lie algebra to the case where the fundamental multiplication is $n-$ary,
$n \geq 2$ (when $n = 2$ the definition agrees with the usual definition of a Lie algebra).
In particular $3-$Lie algebras have close relationships with many important fields in mathematics and
mathematical physics. For example, the study
of supersymmetry, Bagger-Lambert theory, Nambu mechanics or gauge symmetry transformations
of the world-volume theory of multiple coincident M2-branes (see \cite{BL, G, HHM, R, T}). The concept of $n-$Lie superalgebras as generalization of $n-$Lie algebras were initially introduced by Daletskii and Kushnirevich in 1996 (see \cite{DK}). Moreover, Cantarini and Kac gave a more general concept of $n-$Lie superalgebras again in 2010 (see \cite{CK}). $n-$Lie superalgebras are more general structures including $n-$Lie algebras, $n-$ary Nambu-Lie superalgebras, and Lie superalgebras. In 2015, Tao Zhang also introduced  the concept of n-Lie colour algebras, which is the generalization of $n-$Lie superalgebras and  developed  cohomology  theory and deformations of $3-$Lie colour algebras(see \cite{Z}).

Lie-Rinehart algebra was introduced by J. Herz in  \cite{H} and  developed in \cite{Ri} by G. Rinehart, is 
an algebraic generalization of the notion of a Lie algebroid: the space of sections
of a vector bundle is replaced by a module over a ring, a vector field by a derivation
of the ring \cite{Hu}. Some generalizations of Lie-Rinehart algebras, such as Lie-Rinehart superalgebras \cite{Ch} or restricted Lie-Rinehart algebras \cite{D}, have been recently studied.

Determining the inner structure of split algebras  by the techniques of connection of roots becomes more and more meaningful in the area of research in mathematical physics.  It is worth mentioning that these techniques had long been introduced by Calderon, Antonio J, on split Lie algebras with symmetric root systems in \cite{C1}. Recently, in
\cite{ABCS, AC1, CD1, CD2, CNS, CO, CP}, the structure of arbitrary split Lie color algebras, Split regular Hom-Lie algebras, split Leibniz algebras,  split Lie-Rinehart algebras, split $3-$Lie algebras, split $3-$Leibniz
algebras, and arbitrary split Lie algebras of order $3$ 
have been determined by the techniques of connection of roots.

In this paper, we continue to study the structure of  a generalization of  Lie-Rinehart superalgebras to the class which are called  split  $3-$Lie-Rinehart color algebras. The definition, basic structures, actions and crossed modules of $3-$Lie-Rinehart algebras and also  $3-$Lie-Rinehart superalgebras can be fond in \cite{BhCMS, BLW}. A $3-$Lie-Rinehart color algebra is a triple $(\LL, A, \rho),$ where $\LL$ is a $3-$Lie color algebra, $A$ is a commutative, associative color algebra, $\LL$ is an $A-$module, $(A, \rho)$ is a $3-$Lie color algebra $\LL-$module and $\rho(\LL, \LL) \subset Der(A).$ Our goal in this work is to study the inner structure of arbitrary split  $3-$Lie-Rinehart color algebras by the  developing techniques of connections of root systems and weight systems associated to a splitting Cartan subalgebra.  The finding of the present
paper is an improvement and extension of the works conducted in \cite{ABCS}.

We briefly outline the contents of
the paper. In Section 2, we begin by recalling the necessary
background  on split  $3-$Lie-Rinehart color algebr.
Section 3,  develops  techniques of connections of roots for
in the framework of split  $3-$Lie-Rinehart color algebra $(\LL, A, \rho, \e)$ and apply, as a first step, all of these techniques to the study of the inner structure of $\LL.$ In section 4, we get, as a second step, a decomposition of $A$  as the orthogonal
direct sums of adequate ideals. In section 5, we  show that the decompositions of $\LL$ and $A$ as direct sum of ideals, given in Sections 2 and 3 respectively, are closely related. Section 6, devoted to show that  under certain conditions, the given decomposition of $\LL$ and $A$ are by means of the family of their corresponding, graded simple ideals.

Throughout this paper, algebras and vector spaces are over a field $\bbbf$ of characterestic zero, and $A$ denotes an associative and commutative algebra over $\bbbf.$ We also
consider an additive abelian group $G$ with identity zero.

\section{Preliminaries} \setcounter{equation}{0}\

Let us begin with some definitions and results concerning graded algebric structures. For a detailed discussion of this subject, we refer
the reader to the literature \cite{S}.  Let $G$ be any
additive abelian group, a vector space $V$ is
said to be {\em $G-$graded}, if there is a family
$\{V_g\}_{g\in G}$ vector subspaces such that
$V=\bigoplus_{g\in G}V_g.$ An element $v\in V$ is said to
be {\em homogeneous of  degree $g$} if $v\in
V_g,~g\in G,$ and in this case, $g$ is called the {\em
color} of $v.$ As usual, we denote by $|v|$ the color of an element
$v \in V.$ Thus, each homogeneous element $v$ in $V$ determines a
unique group element $|v|\in G$ by $v\in V_{|v|}.$ Fortunately,
we can almost drop the symbol "$|~~|$", since confusion rarely
occurs.

Let $V=\bigoplus_{g\in G}V_g$ and
$W=\bigoplus_{g\in G}W_g$ be two $G-$graded vector
spaces. A linear mapping $f : V\longrightarrow W$ is said to be
{\em homogeneous of degree} $h\in G$ if
$$
f(V_g)\subset W_{g+h}, ~~~\forall g\in G.
$$
If in addition, $f$ is homogeneous of degree zero, namely,
$f(V_g)\subset W_g$ holds for any $g\in G,$ then we call $f$
is  {\em even}.

An algebra $A$ is said to be $G-$graded or color algebra if its underlying
vector space is $G-$graded, i.e.,
$A=\bigoplus_{g\in G}A_g,$ and if $A_g A_h\subset
A_{g+h},$ for $g, h\in G.$ A subalgebra (an ideal) of $A$ is said
to be graded if it is a graded as a subspace of $A.$

Let $B$ be another $G-$graded algebra. A homomorphism $\varphi
: A\longrightarrow B$ of $G-$graded algebras is a homomorphism
of the algebra $A$ into the algebra $B,$ which is an even mapping.

\begin{DEF}\label{bi} Let $G$ be an abelian group. A map $\e :
G\times G\longrightarrow\bbbk\setminus\{0\}$ is called a
skew-symmetric {\em bi-character} on $G$ if for all $g, h, f\in G,$
\begin{itemize}
\item[(i)] $\e(g, h)\e(h, g)=1,$

\item[(ii)] $\e(g+h, f)=\e(g, f) \e(h, f),$

\item[(iii)]  $\e(g, h+f)=\e(g, h)\e(g, f).$
\end{itemize}
\end{DEF}
The definition above implies that in particular, the following
relations hold
$$
\e(g, 0)=1=\e(0, g),~~~~\e(g, g)=1(\hbox{or}~-1),~~~\forall g\in G.
$$
Throughout this paper, if $x$ and $ y$  are homogeneous elements
of a $G-$graded vector space and $|x|$ and $|y|$  which
are in $G$ denote their degrees respectively, then for
convenience, we write $\e(x, y)$ instead of $\e(|x|, |y|).$ It is worth mentioning that,  we
 unless otherwise stated, in the sequel all the graded spaces are over
the same abelian group $G$ and the bi-character will be the
same for all structures.

\begin{DEF}\cite{Z} A $3-$Lie colour algebra consists of a $G-$graded vector space $\LL=\bigoplus_{g\in G}\LL^g$
together with a trilinear map $[., ., .] : \LL\times\LL\times\LL\longrightarrow\LL$ such that the following
conditions are satisfied:
\begin{itemize}
\item[(i)]graded condition:

$|[x_1, x_2, x_3]| = |x_1|+|x_2|+|x_3|;$ 

\item[(ii)] $\e-$skew symmetry:

$[x_1, x_2, x_3]= -\e(x_1, x_2)[x_2, x_1, x_3], [x_1, x_2, x_3]= -\e(x_2, x_3)[x_1, x_3, x_2];$

\item[(iii)] $\e-$fundamental identity:
\begin{eqnarray}\label{0}
[x_1, x_2, [y_1, y_2, y_3]]&=& [[x_1, x_2, y_1], y_2, y_3]\\
\nonumber&+&\e(x_1+x_2, y_1)[y_1, [x_1, x_2, y_2], y_3]\\
\nonumber&+&\e(x_1+x_2, y_1+y_2)[y_1, y_2, [x_1, x_2, y_3]],
\end{eqnarray}
\end{itemize}
for all homogeneous elements $x_1,x_2, y_1,y_2, y_3\in\LL.$
\end{DEF}
Let $(\LL,  [.,., .], \e)$ be a   $3-$Lie  color algebra. A {\em subalgebra} $S$ of  $\LL$ is a $G-$graded subspace $S=\bigoplus_{g\in G} S^g$ of $\LL$ such that $[S, S, S ]\subset S.$ and  A $G-$graded subspace $I=\bigoplus_{g\in G}I^g$
of $\LL$ is called an {\em  ideal} if $[I, \LL, \LL]\subset I.$  A $3-$Lie color algebra $\LL$
is called {\em simple} if its triple product is nonzero and its only ideals are $\{0\}$ and $\LL$. We recall that the
Annihilator of a $3-$Lie color algebra $\LL$ is defined as the set of elements $x$ in $\LL$ such that $[x,\LL,\LL] = 0$. This is an ideal of $\LL$ denoted by $Ann(\LL)$.

\begin{DEF}\cite{BhCMS} Let $(\LL, [., ., .], \e)$ be a $3-$Lie  color algebra  and $V$ be a $G-$graded vector space and $\rho:\LL\times\LL \longrightarrow gl(V)$ be a even linear mapping.
Then $(V,\rho)$ is called a representation of $\LL$ or $V$ is an $\LL-$module  if the following two conditions are satisfied:
\begin{eqnarray}\label{2}
\rho(x_1, x_2)\rho(x_3, x_4) &-&\e(x_1+x_2, x_3+x_4)\rho(x_3, x_4)\rho(x_1, x_2)\\
\nonumber&=& \rho([x_1, x_2, x_3], x_4) -\e(x_3, x_4)\rho([x_1, x_2, x_4], x_3),
\end{eqnarray}
and
\begin{eqnarray}\label{3}
\hspace{7mm}\rho([x_1, x_2, x_3], x_4)&=&  \rho(x_1, x_2)\rho(x_3, x_4) +\e(x_1, x_2+x_3)\rho(x_2, x_3)\rho(x_1, x_4)\\
\nonumber&+&\e(x_1+x_2, x_3)\rho(x_3, x_1)\rho(x_2, x_4).
\end{eqnarray}
\end{DEF}
Next, define
$$
ad : \LL\times\LL\longrightarrow gl(\LL);~~~ad_{x, y}z=[x,y,z] ,~~\forall x, y, z\in\LL.
$$
Tanks to $\e-$fundamental identity, $(\LL, ad)$ is a representation of the  $3-$Lie  color algebra $\LL,$ and it is called the adjoint representation of $\LL.$ One can see that $ad(\LL, \LL)$ is a Lie algebra which is called inner derivation of $\LL.$ We also have by Eq (\ref{2}),
$$
[ad_{x_1, x_2}, ad_{y_1, y_2}]=ad_{[x_1, y_1, x_2], y_2}+ ad_{[x_2, [x_1, y_1, y_2]]}.
$$

\begin{DEF}  Let   $(\LL, [., ., .], \e)$ be a $3-$Lie  color algebra and $A$ be  an asociative commutative color algebra which $(A, \rho)$ is an $\LL-$module. If $\rho(\LL, \LL)\subset Der(A)$ and
\begin{eqnarray}\label{4}
[x, y, a z] =\e(a, x+y) a[x, y, z]+\rho(x, y) a z,~~~\forall x, y, z\in\LL,~~\forall a\in A,
\end{eqnarray}
\begin{eqnarray}\label{5}
\rho(ax, y) =\e(a, x)  \rho(x, ay)= a\rho(x, y),~~~\forall x, y\in\LL,~~\forall a\in A,
\end{eqnarray}	
then $(\LL, A, [., ., .], \rho, \e)$	is called a {\em $3-$Lie-Rinehart color algebra}.
\end{DEF}

\begin{exa}
\begin{itemize}
\item[(1)]  If $\rho=0$ then $(\LL, A, [., ., .], \e)$ is a $3-$Lie $A-$color algebra.
\item[(2)] Let $(\g, [., .], \e)$ be a color Lie algebra and an $A-$module. Let $(A, \rho)$ be a $\g-$module. If $\rho(\g)\subset Der(A)$  and
$$
[x, ay]=\e(a, x)a[x, y]+\rho(x)ay,~~~~\rho(ax)=a\rho(x),~~~~~\forall x, y\in\g,  ~  \forall a\in A,
$$
then $(\g, A, [., .], \rho, \e)$ a Lie-Rinehart color algebra.
\item[(3)] if $G=\bbbz_2$ and $\e(x, y)=(-1)^{x y}$  for any homogeneous elements $x, y\in\LL.$ We then get the notion of $3-$Lie-Rinehart superalgebra (see \cite{BhCMS}).
\end{itemize} 
\end{exa}

\begin{exa} We recall that given a Lie color algebra analogues of color trace
	one can construct a $3-$Lie color algebra. Let $(\LL, [., . ], \e)$ be a Lie color algebra and
	$\tau : \LL\longrightarrow\bbbf$ an even linear form. We say that $\tau$ is a color trace of $\LL$ if $\tau([., .]) = 0.$ For
	any $x_1, x_2, x_3\in\LL,$ we define the $3-$ary bracket by
	\begin{eqnarray}\label{500}
	[x_1, x_2, x_3]_\tau &=& \tau(x_1)[x_2, x_3] - \e(x_1,x_2)\tau(x_2)[x_1, x_3]\\
	\nonumber &+& \e(x_3, x_1+x_2)\tau(x_3)[x_1, x_2]. 
	\end{eqnarray}	
	then $(\LL, [., ., .]_\tau, \e )$ is
	a $3-$Lie colr algebra (see  \cite{A} for super algebras).
	Next, We begin by constructing
	$3-$Lie-Rinehart color algebras starting with a Lie-Rinehart color algebras. Let $(\LL, A, [., .], \rho, \e)$ be a Lie-Rinehart color algebra and $\tau$ is a color trace. If the condition
	$$
	\tau(ax)y = \tau(x)ay,
	$$ 
	is satisfied for any  $x, y\in\LL, a\in A$, then $(\LL, A, [., ., .]_\tau, \rho_\tau)$ is a $3-$Lie-Rinehart color algebra, where $[., ., .]_\tau$ is defined as Eq.(\ref{500}) and $\rho_\tau$ is defined by
	$$
	\rho_\tau :\LL\times\LL\longrightarrow gl(V) ;~~\rho_\tau(x, y)=\tau(x)\rho(y)-\e(x, y)\tau(y)\rho(x),~~\forall x, y\in\LL.
	$$
	(see Theorem 2.1 in \cite{BhCMS} for super algebras)	
\end{exa}

\begin{DEF}\label{ideal} Let $(\LL, A, [., ., .], \rho, \e)$	be a $3-$Lie-Rinehart color algebra.
\begin{itemize}
\item[(1)] If $S$ is a $3-$Lie color subalgebra  of  $\LL$ satisfying $A S\subset S,$ then $$(S, A, [., ., .]|_{S\times S}, \rho|_{S\times S}, \e)$$ is a  $3-$Lie-Rinehart color algebra which is called a {\em subalgebra} of the   $3-$Lie-Rinehart color algebra  $(\LL, A, [., ., .], \rho, \e).$

\item[(2)]  If $I$ is a $3-$Lie color ideal  of  $\LL$ satisfying $A I\subset I$ and $\rho(I, I)(A)\subset  I,$ then $$(I, A, [., ., .]|_{I\times I}, \rho|_{I\times I}, \e)$$ is a  $3-$Lie-Rinehart color algebra which is called an {\em ideal} of the   $3-$Lie-Rinehart color algebra  $(\LL, A, [., ., .], \rho, \e).$ As in \cite{BhCMS}, Proposition 2. 4, $\ker\rho$ is an ideal of $(\LL, A, [., ., .], \rho, \e).$ 
\item[(3)] If a  $3-$Lie-Rinehart color algebra  $(\LL, A, [., ., .], \rho, \e)$ cannot be decomposed into the direct sum of two non-zero ideals, then it is called an {\em indecomposabel}  $3-$Lie-Rinehart color algebra. We also  say that  $(\LL, A, [., ., .], \rho, \e)$ is {\em simple} if $[\LL, \LL, \LL]\neq 0,~AA\neq 0, \LL A\neq 0$ and its only ideals are $\{0\}$, $\LL$ and $\ker \rho.$ 
\end{itemize}
\end{DEF}
Denotes;
$$
Z_\LL(A):=\{a\in A : ax=0~\forall x\in\LL\},~~~~~Z_\rho(\LL):=\{x\in\LL : [x, \LL, \LL]=0,~ \rho(x, \LL)=0\}.
$$
Note that $Z_\rho(\LL)=Ann(\LL)\cap\ker\rho.$

From now on,   $(\LL, A, [., ., .], \rho, \e)$ denotes a  $3-$Lie-Rinehart color algebra. We introduce the class of split algebras in the framework of $3-$Lie-Rinehart color algebras as in \cite{ABCS}.

\begin{DEF} Let  $(\LL, A, [., ., .], \rho, \e)$ be a $3-$Lie-Rinehart color algebra. If there
exist a maximal abelian subalgebra $\hh$ of $\LL$ satisfying that 
\begin{equation}\label{30}
\LL=\hh\oplus(\bigoplus_{\a\in\Pi}\LL_\a),~~\hbox{where}~~ \pi = \{0\neq\a\in(\hh\times\hh)^* ~ |~ \LL_\a\neq 0\},
\end{equation}
\begin{equation}\label{31}
\LL_\a = \{x \in \LL~ |~ [h_1, h_2, x] = \a(h_1, h_2)(x), \forall h_1, h_2 \in\hh\}, \LL_0 = \hh.
\end{equation}
And
\begin{equation}\label{32}
A=A_0\oplus(\bigoplus_{\lam\in\Lam}A_\lam),~~\hbox{where}~~ \Lam = \{0\neq\lam\in(\hh\times\hh)^* ~ |~ A_\lam\neq 0\},
\end{equation}
\begin{equation}\label{33}
A_\lam= \{a \in A ~ |~ \rho(h_1, h_2)(a) = \lam(h_1, h_2)(a), \forall h_1, h_2 \in\hh \},
\end{equation}

\begin{equation*}
A_0=\{a \in A ~ |~ \rho(\hh, \hh)(a) =0 \}.
\end{equation*}
then  $(\LL, A, [., ., .], \rho, \e)$ is called a {\em split}  $3-$Lie-Rinehart color algebra, and $\hh$ is a {\em splitting Cartan subalgebra} of $\LL$, $\Pi$ and $\Lam$ are called the {\em root system} of $\LL$, and
$\Lam$ the {\em weight system} of $A$ associated to $\hh$, respectively. Linear subspaces $\LL_\a (\a\in\Pi)$ and $A_\lam (\lam\in\Lam)$ are called the {\em root space} of $\LL$ and the {\em weight space} of $A,$ respectively.
\end{DEF}

For convenience, in the rest of the paper, the  split  $3-$Lie-Rinehart color algebra   $(\LL, A, [., ., .], \rho, \e)$ is simply denoted by $(\LL, A),$ and denote
$$
-\pi = \{-\a~ |~ \a\in\pi\},~ -\Lam = \{-\lam ~ |~  \lam\in\Lam\},~~~~ \pm\pi = \pi\cup -\pi,~ \pm\Lam = \Lam\cap-\Lam.
$$

We recall some properties of split   $3-$Lie-Rinehart color algebra that can be found in \cite{Kh1} for a split regular involutive Hom-Lie color algebras. 

\begin{lem}\label{2.6} Let
$(\LL=\bigoplus_{g\in G}\LL^g, A=\bigoplus_{g\in G}A_g)$ be a split 
$3-$Lie-Rinehart color algebra, with root space decomposition
$\LL=\hh\oplus(\bigoplus_{\a\in\Pi}\LL_\a)$ and weight space decomposition $A=A_0\oplus(\bigoplus_{\lam\in\Lam}A_\lam).$ Then
\begin{itemize}

\item[(1)] for any $\a\in\Pi\cup\{0\},$ we have
$\LL_\a=\bigoplus_{g\in G}\LL_\a^g,$ where
$\LL_\a^g=\LL^g\cap\LL_\a,$ and $\LL_0=\hh,$
\item[(2)] $\hh^g=\LL_0^g.$ In particular, $\hh^0=\LL_0^0,$
\item[(3)] $\LL^0$ is a split $3-$Lie color algebra, with respect to
$\hh^0=\LL_0^0,$ with root space decomposition
$\LL^0=\hh^0\oplus(\bigoplus_{\a\in\Pi}\LL^0_\a).$
\item[(4)] for any $\lam\in\Lam\cup\{0\},$ we have
$A_\lam=\bigoplus_{g\in G}A_\lam^g,$ where
$A_\lam^g=A^g\cap A_\lam.$ We also have $A^g=A_0^g,$ in particular $A^0=A_0^0.$
\end{itemize}
\end{lem}
\noindent {\bf Proof.} Similar to Lemma 2.6 in \cite{Kh1}.\qed

If $(\LL, A)$ is a split   $3-$Lie-Rinehart color algebra,  with root space decomposition
$\LL=\hh\oplus(\bigoplus_{\a\in\Pi}\LL_\a)$ and weight space decomposition $A=A_0\oplus(\bigoplus_{\lam\in\Lam}A_\lam),$ taking into account 
Lemma \ref{2.6}, we then write
\begin{equation}\label{00}
\LL=\bigoplus_{g\in G}(\hh^g\oplus(\bigoplus_{\a\in\Pi}\LL^g_\a))=\hh^0\oplus(\bigoplus_{g\in G}\bigoplus_{\a\in\Pi^g}\LL^g_\a),
\end{equation}
\begin{equation}\label{00.}
A=\bigoplus_{g\in G}(A^g\oplus(\bigoplus_{\lam\in\Lam}A^g_\lam))=A^0\oplus(\bigoplus_{g\in G}\bigoplus_{\lam\in\Lam^g}A^g_\lam),
\end{equation}

we denote by $\Pi^g : =\{\a\in\Pi ~|~\LL^g_\a\neq 0\}$  and  $\Lam^g : =\{\lam\in\Lam ~|~A^g_\a\neq 0\},$ for any $g\in G.$ Then $\Pi=\bigcup_{g\in G}\Pi^g$ and $\Lam=\bigcup_{g\in G}\Lam^g.$

\begin{pro}\label{2.7} For any $\a_1, \a_2, \a_3\in\Pi\cup\{0\},~ \lam_1, \lam_2\in\Lam\cup\{0\} $ and any $g_1, g_2, g_3\in G,$ the
following assertions hold.
\begin{itemize}
\item[(1)] If  $[\LL^{g_1}_{\a_1}, \LL^{g_2}_{\a_2}, \LL^{g_3}_{\a_3}]\neq 0,$ then $\a_1+\a_2+ \a_3\in\Pi^{g_1+g_2+g_3}\cup\{0\}$ and  $$[\LL^{g_1}_{\a_1}, \LL^{g_2}_{\a_2}, \LL^{g_3}_{\a_3}]\subset\LL^{g_1+g_2+g_3}_{\a_1+\a_2+\a_3}.$$

\item[(2)] If $A_{\lam_1}^{g_1}A_{\lam_2}^{g_2}\neq 0,$ then $\lam_1+ \lam_2\in\Lam^{g_1+g_2}\cup\{0\}$ and $A_{\lam_1}^{g_1}A_{\lam_2}^{g_2}\subset A_{\lam_1+ \lam_2}^{g_1+g_2}.$

\item[(3)] If $A_{\lam_1}^{g_1}\LL^{g_2}_{\a_2}\neq 0,$ then $\lam_1+ \a_2\in\Pi^{g_1+g_2}\cup\{0\}$ and  $A_{\lam_1}^{g_1}\LL^{g_2}_{\a_2}\subset \LL^{g_1+g_2}_{\lam_1+\a_2}.$

\item[(4)] If $\rho(\LL^{g_1}_{\a_1}, \LL^{g_2}_{\a_2})(A_{\lam_3}^{g_3})\neq 0,$ then $\a_1+\a_2+\lam_3\in\Lam^{g_1+g_2+g_3}\cup\{0\}$ and $$\rho(\LL^{g_1}_{\a_1}, \LL^{g_2}_{\a_2})(A_{\lam_3}^{g_3})\subset A_{\a_1+\a_2+\lam_3}^{g_1+ g_2+ g_3}.$$
\end{itemize}
\end{pro} 
\noindent {\bf Proof.} 

(1) For any $h_1, h_2\in\hh^0,~~y_i\in\LL_{\a_i}^{g_i},~i=1, 2, 3$ from  $\e-$fundamental identity, we have
\begin{eqnarray*}
[h_1, h_2, [y_1, y_2, y_3]]&=& [[h_1, h_2, y_1], y_2, y_3]+\e(0, y_1)[y_1, [h_1, h_2, y_2], y_3]\\
&+&\e(0, y_1+y_2)[y_1, y_2, [h_1, h_2, y_3]]\\
&=&[\a_1(h_1, h_2)y_1, y_2, y_3]+\e(0, y_1)[y_1, \a_2(h_1, h_2) y_2, y_3]\\
&+&\e(0, y_1+y_2)[y_1, y_2, \a_3(h_1, h_2) y_3]\\
&=&(\a_1+\a_2+\a_3)(h_1, h_2) [y_1, y_2, y_3].
\end{eqnarray*}
Therefore, we get  $\a_1+\a_2+ \a_3\in\Pi^{g_1+g_2+g_3}g\cup\{0\}$ and  $[\LL^{g_1}_{\a_1}, \LL^{g_2}_{\a_2}, \LL^{g_3}_{\a_3}]\subset\LL^{g_1+g_2+g_3}_{\a_1+\a_2+\a_3}.$

(2)  For any $a_1\in A_{\lam_1}^{g_1}, a_2\in A_{\lam_2}^{g_2},$ since $\rho(h_1, h_2)\in Der(A), ~\forall h_1, h_2\in\hh^0$ we have
\begin{eqnarray*}
\rho(h_1, h_2) (a_1 a_2)&=& \rho(h_1, h_2) (a_1) a_2+a_1 \rho(h_1, h_2) (a_2)\\
&=& \lam_1(h_1, h_2) (a_1) a_2+a_1 \lam_2(h_1, h_2) (a_2)\\
&=&(\lam_1+\lam_2)(h_1, h_2) (a_1 a_2).
\end{eqnarray*}
Thus $\lam_1+ \lam_2\in\Lam^{g_1+g_2}\cup\{0\}$ and $a_1 a_2\in A_{\lam_1+ \lam_2}^{g_1+g_2}.$

(3)  For any $a\in A_{\lam_1}^{g_1},~x\in\LL^{g_2}_{\a_2}$ and $h_1, h_2\in\hh^0,$ by Eq. (\ref{4}) we have
\begin{eqnarray*}
[h_1, h_2, ax ] &=& \e(a, 0)a[h_1, h_2, x]+\rho(h_1, h_2) a x\\
&=&a\a_2(h_1, h_2)x+\lam_1(h_1, h_2)a x\\
&=&(\a_2+\lam_1)(h_1, h_2)(ax).
\end{eqnarray*}
Therefore, we get $\lam_1+ \a_2\in\pi^{g_1+g_2}\cup\{0\}$ and $ax\in\LL^{g_1+g_2}_{\lam_1+\a_2}.$

(4) Since $A$ is an $\LL-$module, for any $a\in A_{\lam_3}^{g_3},~x_i\in\LL^{g_i}_{\a_i}(i=1, 2)$ and   $h_1, h_2\in\hh^0,$ we have

\begin{eqnarray*}
\rho(h_1, h_2)(\rho(x_1, x_2)a) &=& (\rho(h_1, h_2)\rho(x_1, x_2))(a)\\
&=&\e(0, x_1+x_2)\rho(x_1, x_2)\rho(h_1, h_2)(a)+\rho([h_1, h_2, x_1], x_2)a\\
&-&\e(x_1, x_2)\rho([h_1, h_2, x_2], x_1)a\\
&=&\rho(x_1, x_2)\lam_1(h_1, h_2)(a)+\rho(\a_1(h_1, h_2)x_1, x_2)a\\
&-&\e(x_1, x_2)\rho(\a_2(h_1, h_2)x_2, x_1)a\\
&=&(a_1+\a_2+\lam_1)(h_1, h_2)\rho(x_1, x_2)(a).
\end{eqnarray*}
Thus $\a_1+\a_2+\lam_3\in\Lam^{g_1+g_2+g_3}\cup\{0\}$ and $\rho(x_1, x_2)(a)\in A_{\a_1+\a_2+\lam_3}^{g_1+ g_2+ g_3}.$\qed

\section{Connections of roots and decompositions} \setcounter{equation}{0}\
In this section, we begin by developing the techniques
of connections of roots in the setting as  \cite{ABCS}. Let  $(\LL, A)$ be a split   $3-$Lie-Rinehart color algebra,  with root space decomposition
$\LL=\hh^0\oplus(\bigoplus_{g\in G}\bigoplus_{\a\in\Pi^g}\LL^g_\a),$ and  with  symmetric root system $\Pi=\bigcup_{g\in G}\Pi^g.$ 

\begin{DEF}\label{conn}
Let $\a, ~\b$ be two non zero roots in $\Pi.$ We say that {\em $\a$
is connected to $\b$} and denoted by $\a\sim\b$ if there exists a
family
$$
\{\a_1, \a_2, \a_3, ..., \a_{2k+1}\}\subset\pm\Pi\cup\pm\Lam\cup\{0\},
$$
satisfying the following conditions;\\

If $k=1:$
\begin{itemize}
\item[(1)] $\b=\pm\a_1.$
\end{itemize}

If $k\geq2:$
\begin{itemize}
\item[(1)] $\a_1+(\a_2+\a_3)\in\pm\pi,\\
\a_1+ (\a_2+( \a_3+(\a_4+\a_5)))\in\pm\Pi,\\
...\\
\a_i+ \sum_{j=1}^i(\a_{2j}+ \a_{2j+1})\in\pm\Pi,~~0< i< k.$

\item[(2)] $\a_1+ \sum_{j=1}^k(\a_{2j}+ \a_{2j+1})\in\{\b, -\b\}.$
\end{itemize}
The family $\{\a_1, \a_2, \a_3, ..., \a_{2k+1}\}$ is called a
{\em connection} from $\a$ to $\b.$
\end{DEF}
\begin{pro}\label{rela}
The relation $\sim$ in $\Pi$ defined by
$$
\a\sim\b~\hbox{ if and only if}~ \a~\hbox{ is connected to}~ \b,
$$
is an equivalence relation.
\end{pro}
\noindent {\bf Proof.} The proof is vertically identical to the
proof  of Proposition 4.5 in \cite{BLW}. \qed

\begin{lem}\label{8.} For any $\gamma, \mu\in\pm\Pi\cup\pm\Lam\cup\{0\},$ if $\a\sim\b$ and $\a+\gamma+\mu\in\pi,$ then $\b\sim\a+\gamma+\mu.$ 
\end{lem}
\noindent {\bf Proof.} Considering the connection $\{\b, \gamma, \mu\}$ we get  $\a\sim\a+\gamma+\mu.$ Now, by transitivity $\b\sim\a+\gamma+\mu.$\qed

By the  Proposition \ref{rela}, we can consider the equivalence relation
in $\Pi$ by the connection relation $\sim$ in $\Pi.$ So we denote
by
$$
\Pi/\sim :=\{[\a] : \a\in\Pi\},
$$
where $[\a]$ denotes  the set of non zero roots of $\LL$ which are
connected to $\a.$ 

Our next goal in this section is to associate an adequate ideal
$I_{[\a]}$ of $\LL$ to any $[\a].$ For a fixed  $\a\in\Pi,$ we define
$$
I_{0, [\a]} :=(\sum_{\substack{\b\in[\a]\\ -\b\in\Lam}} A_{-\b}\LL_\b)+(\sum_{\substack{\b, \gamma, \mu\in[\a]\\ \b+\gamma+\mu=0}}[\LL_\b,
\LL_{\gamma}, \LL_\mu]).
$$
In order to graded case, 
\begin{eqnarray}\label{119.}
I_{0, [\a]} : =(\sum_{\substack{\b\in [\a]\\ -\b\in\Lam\\ g, g'\in G}} A_{-\b}^g\LL_\b^{g'})
+(\sum_{\substack{\b, \gamma, \mu\in [\a]\\ \b+\gamma+\mu=0\\ g, g', g''\in G}}[\LL_\b^g,
\LL_{\gamma}^{g'}, \LL_\mu^{g''}]).
\end{eqnarray}
Next, we define
\begin{equation}\label{38'}
\vv_{[\a]}
:=\bigoplus_{\b\in[\a]}\LL_\b=\bigoplus_{g\in G}\bigoplus_{\b\in[\a]}\LL^g_\b.
\end{equation}
Thanks to Proposition \ref{2.7}, $I_{0, [\a]}\subset\hh$ and $I_{0, [\a]}\cap\vv_{[\a]}=0.$
Finally, we denote by $I_{[\a]}$ the direct sum of the two graded
subspaces above, that is,
\begin{equation}\label{38''}
I_{[\a]} :=I_{0, [\a]}\oplus\vv_{[\a]}.
\end{equation}

\begin{pro}\label{subalg} For any $[\a]\in\Pi/\sim$
 the following assertions hold.
 \begin{itemize}
 \item[(1)] $[I_{[\a]}, I_{[\a]}, I_{[\a]}]\subset I_{[\a]}.$
 \item[(2)]$A I_{[\a]}\subset I_{[\a]} $
 \item[(3)]$\rho( I_{[\a]}, I_{[\a]})(A)\LL\subset I_{[\a]}$		
 \end{itemize}
\end{pro}
\noindent {\bf Proof.} (1) By
the fact $\LL_0=\hh$ and Eq. (\ref{119.}), it is clear that $[I_{0,
[\a]},  I_{0, [\a]},  I_{0, [\a]}]\subset[\hh, \hh, \hh]=0,$ and we can write
\begin{eqnarray}\label{200}
\nonumber[I_{[\a]}, I_{[\a]}, I_{[\a]}]&=&[I_{0, [\a]}\oplus\vv_{[\a]},
I_{0, [\a]}\oplus\vv_{[\a]}, I_{0, [\a]}\oplus\vv_{[\a]}]\\
&\subset&[I_{0, [\a]}, \vv_{[\a]}, I_{0, [\a]}]+[I_{0, [\a]}, \vv_{[\a]}, \vv_{[\a]}]+[\vv_{[\a]}, \vv_{[\a]}, \vv_{[\a]}].
\end{eqnarray}
Let us consider the first summand in (\ref{200}). Since $ I_{0, [\a]}\subset\hh,$ for $\b\in[\a]$ and $g\in G,$ by Proposition \ref{2.7}, one gets  $[I_{0,
[\a]}, I_{0, [\a]}, \LL^g_\b]\subset\LL^g_\b.$  Hence,
\begin{equation}\label{38}
[I_{0, [\a]}, \vv_{[\a]}, I_{0, [\a]}]\subset \vv_{[\a]}.
\end{equation}
Consider now the second  summand in (\ref{200}). Given $\b, \gamma\in [\a]$ and $g, g'\in G$ such that $[I_{0, [\a]}, \LL_\b^g, \LL_\gamma^{g'}]\neq 0,$ then  $[I_{0, [\a]}, \LL_\b^g, \LL_\gamma^{g'}]\subset\LL_{\b+\gamma}^{g+g'},$  by Proposition \ref{2.7}. If $\b+\gamma\in\Pi^{g+g'},$ and $\{\b, \gamma, 0\}$ is a connection from $\b$ to $\b+\gamma,$ by Lemma \ref{8.}. We have
$$
[I_{0, [\a]}, \LL_\b^g, \LL_\gamma^{g'}]\subset\vv_{[\a]}\subset I_{[\a]}.
$$
Therefor,
\begin{equation}\label{39}
[I_{0, [\a]}, \vv_{[\a]}, \vv_{[\a]}]\subset I_{[\a]}.
\end{equation}
Finally, consider the last summand in (\ref{200}). Let  $\b, \gamma, \mu\in [\a]$ and $g, g', g''\in G$ such that $[\LL_\b^g, \LL_\gamma^{g'}, \LL_\mu^{g''}]\neq 0,$ and by Proposition \ref{2.7} we have  $[\LL_\b^g, \LL_\gamma^{g'}, \LL_\mu^{g''}]\subset\LL_{\b+\gamma+\mu}^{g+g'+g''}.$ 
If $\b+\gamma+\mu=0,$  we get
$$
[\LL^g_\b, \LL^{g'}_\gamma, \LL_\mu^{g''}]\subset I_{0, [\a]}\subset I_{[\a]}. 
$$
Suppose  $0\neq\b+\gamma+\mu,$ by Lemma \ref{8.}, one gets
$\b+\gamma+\mu\in\Pi^{g+g'+g''},$ therefore, $\{\b, \gamma, \mu\}$ is
a connection from $\b$ to $\b+\gamma+\mu.$ The
transitivity of $\sim$ gives us
$\b+\gamma+\mu\in[\a]$ and so $[\LL^g_\b,
\LL^{g'}_\gamma, \LL_\mu^{g''}]\subset\LL^{g+g'+g''}_{\b+\gamma+\mu}\subset\vv_{[\a]}.$
Hence,
\begin{equation}\label{40}
[\vv_{[\a]}, \vv_{[\a]}, \vv_{[\a]}]\subset I_{[\a]}.
\end{equation}
From Eqs. (\ref{38}), (\ref{39}),  and (\ref{40}), we conclude that
\begin{equation}\label{26}
[I_{[\a]}, I_{[\a]}, I_{[\a]}]\subset I_{[\a]}.
\end{equation}

(2) Observe that
\begin{eqnarray*}
A I_{[\a]}=(A_0\oplus(\bigoplus_{\substack{\lam\in\Lam\\ g\in G}}A_\lam^g))((\sum_{\substack{\b\in [\a]\\ -\b\in\Lam\\ g, g'\in G}} A_{-\b}^g\LL_\b^{g'}
+\sum_{\substack{\b, \gamma, \mu\in [\a]\\ \b+\gamma+\mu=0\\ g, g', g''\in G}}[\LL_\b^g,
\LL_{\gamma}^{g'}, \LL_\mu^{g''}])
\oplus\bigoplus_{\substack{\b\in[\a]\\ g\in G}}\LL_\b^g)
\end{eqnarray*}
We discuss it in six cases:

{\bf Case 1.} For any $\b\in [\a],~g, g'\in G,$ if $-\b\in\Lam^h$ for some $h\in G$ then by Proposition \ref{2.7} and the fact that $\LL$ is an $A-$module, we have
$$
A_0(A_{-\b}^g\LL_\b^{g'})=(A_0 A_{-\b}^g)\LL_\b^{g'}\subset A_{-\b}^g\LL_\b^{g'}\subset I_{0, [\a]}.
$$
Therefore,
$$
A_0(A_{-\b}^g\LL_\b^{g'})\subset I_{[\a]}.
$$
{\bf Case 2.} For any $\b, \gamma, \mu\in [\a],$ with $\b+\gamma+\mu=0,$ and $g, g', g''\in G$ Acording to Eq. (\ref{4}), for any $a_0\in A_0,~ x\in \LL_\b^g,~y\in\LL_\gamma^{g'}$ and $z\in\LL_\mu^{g''},$ we have 
\begin{eqnarray*}
\e(a_0, x+y)a_0[x, y, z]&=&[x, y, a_0 z]-\rho(x, y) a_0 z\\
&\in&[\LL_\b^g, \LL_\gamma^{g'}, A_0 \LL_\mu^{g''}]+\rho(\LL_\b^g, \LL_\gamma^{g'})(A_0)\LL_\mu^{g''}.
\end{eqnarray*}
If $A_{-\mu}^{g+g'}\neq 0$ (otherwise is trivial), then  Proposition \ref{2.7}, gives us $-\mu\in\Lam{g+g'},$ and 
$$
\e(a_0, x+y)a_0[x, y, z]\in[\LL_\b^g, \LL_\gamma^{g'}, \LL_\mu^{g''}]+A_{-\mu}^{g+g'}\LL_\mu^{g''}\subset  I_{0, [\a]}.
$$ 
Therefore,
$$
A_0[\LL_\b^g, \LL_\gamma^{g'}, \LL_\mu^{g''}]\subset I_{[\a]}.
$$

{\bf Case 3.} For any $\b\in [\a],~g\in G,$  by  Proposition \ref{2.7}, we have  
$$
A_0\LL_\b^g\subset\LL_\b^g\subset \vv_{[\a]}.
$$
Hence
$$
A_0\LL_\b^g\subset\LL_\b^g\subset I_{[\a]}.
$$
{\bf Case 4.} For $\lam\in\Lam,~\b\in [\a],$ and $g, g', g''\in G,$ if $-\b\in\Lam^h$ for some $h\in G$ then
$$
A_\lam^{g'}(A_{-\b}^{g''} \LL_\b^g)=(A_\lam^{g'} A_{-\b}^{g''}) \LL_\b^g\subset A_{\lam-\b}^{g'+g''}\LL_\b^g.
$$
If $\lm-\b\in\Lam^h$ for some $h\in G$ and $\LL_\lam^{g+g'+g''}\neq 0$ (otherwise is trivial),  by Proposition \ref{2.7}, $\{\b, \lam-\b, 0\}$ is a connection from $\b$ to $\lam$ and so $\lam\in [\a].$ Hence,
$$
A_\lam^{g'}(A_{-\b}^{g''} \LL_\b^g)\subset\LL_\lam^{g+g'+g''}\subset \vv_{[\a]}.
$$
Therefore,
$$
A_\lam^{g'}(A_{-\b}^{g''} \LL_\b^g)\subset I_{[\a]}.
$$
{\bf Case 5.}  For any $\b, \gamma, \mu\in [\a]$ with $\b+\gamma+\mu=0,$ and  $g_1, g_2, g_3\in G,$ if $\lam\in\Lam,~g\in G.$  Acording to Eq. (\ref{4}), for any $a\in A_\lam^g,~ x\in \LL_\b^{g_1},~y\in\LL_\gamma^{g_2}$ and $z\in\LL_\mu^{g_3},$ we have
\begin{eqnarray*}
\e(a, x+y)a[x, y, z]&=&[x, y, a z]-\rho(x, y) a z\\
&\in&[\LL_\b^{g_1}, \LL_\gamma^{g_2}, A_\lam^g \LL_\mu^{g_3}]+\rho(\LL_\b^{g_1}, \LL_\gamma^{g_2})(A_\lam^g)\LL_\mu^{g_3}\\
&\subset&[\LL_\b^{g_1}, \LL_\gamma^{g_2}, \LL_{\lam+\mu}^{g+g_3}]+A_{\lam-\mu}^{g+g_1+g_2}\LL_\mu^{g_3}.
\end{eqnarray*}
If $\LL_\lam^{g+g_1+g_2+g_3}\neq 0$ and $ \LL_{\lam+\mu}^{g+g_3}\neq 0$ (otherwise is trivial), then $\lam\in\Pi^{g+g_1+g_2+g_3}$ and $\lam+\mu\in\Pi^{g+g_3}.$ By  Proposition \ref{2.7}, we get $\lam\sim\mu,$ so $\lam\in [\a].$ Therefore,
$$
A_\lam^g[\LL_\b^{g_1}, \LL_\gamma^{g_2}, \LL_\mu^{g_3}]\subset I_{[\a]}.
$$

{\bf Case 6.} For any $\lam\in\Lam,~\b\in [\a]$ and $g, g'\in G,$ By  Proposition \ref{2.7}, we get 
$$
A_\lam^g\LL_\b^{g'}\subset \LL_{\lam+\b}^{g+g'}\subset \vv_{[\a]}.
$$
Hence,
$$
A_\lam^g\LL_\b^{g'}\subset I_{[\a]}.
$$
Now, summarizing a discussion of above six cases, we get the result.

(3) By part (2) and Eq. (\ref{4}) we have

$$
\rho( I_{[\a]}, I_{[\a]})(A)\LL\subset [I_{[\a]}, I_{[\a]}, A\LL]+A[I_{[\a]}, I_{[\a]}, \LL]\subset I_{[\a]}.\qed
$$

\begin{pro}\label{12345} Let $[\a], [\b], [\gamma]\in\Pi/\sim$ be different from each other. Then 
$$
[I_{[\a]}, I_{[\b]}, I_{[\gamma]}]=0,~~~\hbox{and~~}~~~~[I_{[\a]}, I_{[\a]}, I_{[\b]}]=0.
$$
\end{pro}
\noindent {\bf Proof.} Since $[I_{0, [\a]}, I_{0, [\b]}, I_{0, [\gamma]}]\subset[\hh, \hh, \hh]=0,$ we have
\begin{eqnarray}\label{205}
[I_{[\a]}, I_{[\b]}, I_{[\gamma]}]&=&[I_{0, [\a]}\oplus\vv_{[\a]}, I_{0, [\b]}\oplus\vv_{[\b]}, I_{0, [\gamma]}\oplus\vv_{[\gamma]}]\\
\nonumber&\subset&[I_{0, [\a]}, I_{0, [\b]}, \vv_{[\gamma]}]+[I_{0, [\a]}, \vv_{[\b]}, I_{0,
[\gamma]}]\\
\nonumber&+&[I_{0, [\a]}, \vv_{[\b]}, \vv_{[\gamma]}]+[\vv_{[\a]}, I_{0, [\b]},  I_{0, [\gamma]}]\\
\nonumber&+&[\vv_{[\a]}, I_{0, [\b]}, \vv_{[\gamma]}]+[\vv_{[\a]}, \vv_{[\b]},  I_{0, [\gamma]}]\\
\nonumber&+&[\vv_{[\a]}, \vv_{[\b]}, \vv_{[\gamma]}]. 
\end{eqnarray}
First, we consider the item $[\vv_{[\a]}, \vv_{[\b]}, \vv_{[\gamma]}]$ in Eq. (\ref{205}).  Suppose that there exist $\a_1\in [\a], \b_1\in[\b], \gamma_1\in[\gamma]$ and $g_1, g_2, g_3\in G$ such that
$[\LL_{\a_1}^{g_1}, \LL_{\b_1}^{g_2}, \LL_{\gamma_1}^{g_3}]\neq 0.$  By  Proposition \ref{2.7}, $\a_1+\b_1+\gamma_1\in\Pi^{g_1+g_2+g_3}.$  Since $\a\sim\a_1,~~\a_1+\b_1+\gamma_1\in\Pi^{g_1+g_2+g_3},$ and Lemma \ref{8.}  give us $\a\sim\a_1+\b_1+\gamma_1.$ Similarly $\b\sim\a_1+\b_1+\gamma_1.$ By transitivity, $\a\sim\b,$ which is a contradiction. Hence, 
\begin{equation}\label{149}
[\vv_{[\a]}, \vv_{[\b]}, \vv_{[\gamma]}]=\{0\}.
\end{equation}
Next, we consider the item $[I_{0, [\a]}, \vv_{[\b]}, \vv_{[\gamma]}]$ in Eq. (\ref{205}). We have
\begin{eqnarray*}
[I_{0, [\a]}, \vv_{[\b]}, \vv_{[\gamma]}] =(\sum_{\substack{\b\in [\a]\\ -\b\in\Lam\\ g, g'\in G}} A_{-\b}^g\LL_\b^{g'})
+(\sum_{\substack{\b, \gamma, \mu\in [\a]\\ \b+\gamma+\mu=0\\ g, g', g''\in G}}[\LL_\b^g,
\LL_{\gamma}^{g'}, \LL_\mu^{g''}]), \vv_{[\b]}, \vv_{[\gamma]}].
\end{eqnarray*}
For $\b_1\in[\b], \gamma_1\in[\gamma],$ and $x_i\in\LL_{\a_i}^{g_i} (i=1, 2, 3),~ y_2\in\LL_\b^{h},~ y_3\in\LL_{\gamma}^k$ by $\e-$fundamental identity,
\begin{eqnarray*}
 [[x_1, x_2, x_3], y_2, y_3]]&=& [[x_1,  y_2, y_3], x_2, x_3]
 +\e(x_1+x_2, y_2)[x_1, [x_2, y_2, y_3], x_3]\\
 \nonumber&+&\e(x_1+x_2, y_2+y_3)[x_1, x_2, [x_3, y_2, y_3]]\\
 &\subset& [[\LL_{\a_1}^{g_1},  \LL_\b^h, \LL_\gamma^k], \LL_{\a_2}^{g_2}, \LL_{\a_3}^{g_3}]
 +[\LL_{\a_1}^{g_1}, [\LL_{\a_2}^{g_2}, \LL_\b^h, \LL_\gamma^k]], \LL_{\a_3}^{g_3}]\\
 &+&[\LL_{\a_1}^{g_1}, \LL_{\a_2}^{g_2}, [\LL_{\a_3}^{g_3}, \LL_\b^h, \LL_\gamma^k]].
\end{eqnarray*}
Using Eq.(\ref{149}), we get
$$
[\LL_{\a_1}^{g_1},  \LL_\b^h, \LL_\gamma^k]= [\LL_{\a_2}^{g_2}, \LL_\b^h, \LL_\gamma^k]= [\LL_{\a_3}^{g_3}, \LL_\b^h, \LL_\gamma^k]=0.
$$
Therefore,
$$
[(\sum_{\substack{\b, \gamma, \mu\in [\a]\\ \b+\gamma+\mu=0\\ g, g', g''\in G}}[\LL_\b^g,
\LL_{\gamma}^{g'}, \LL_\mu^{g''}]), \vv_{[\b]}, \vv_{[\gamma]}]=0.
$$
Now, if there exist $\b_1\in[\b], \gamma_1\in[\gamma]$ such that $[A_{-\a_1}^{g}\LL_{\a_1}^{g'},\LL_{\b_1}^h, \LL_{\gamma_1}^k]\neq 0.$ By Eq. (\ref{4}), we have

\begin{eqnarray*}
[A_{-\a_1}^{g}\LL_{\a_1}^{g'},\LL_{\b_1}^h, \LL_{\gamma_1}^k]&=&[\LL_{\b_1}^h, \LL_{\gamma_1}^k, A_{-\a_1}^{g}\LL_{\a_1}^{g'}]\\
&\subset&\ A_{-\a_1}^g[\LL_{\a_1}^{g'},\LL_{\b_1}^h, \LL_{\gamma_1}^k]+\rho(\LL_{\b_1}^h, \LL_{\gamma_1}^k)(A_{-\a_1}^g)\LL_{\a_1}^{g'}.
\end{eqnarray*}
By Eq.(\ref{149}), we have $[\LL_{\a_1}^{g'},\LL_{\b_1}^h, \LL_{\gamma_1}^k]=0,$ and   Proposition \ref{2.7}, gives us

$$
\rho(\LL_{\b_1}^h, \LL_{\gamma_1}^k)(A_{-\a_1}^g)\LL_{\a_1}^{g'}\subset A_{{\b_1}+{\gamma_1}-\a_1}^{g+h+k}\LL_{\a_1}^{g'}\neq 0,
$$
then $A_{{\b_1}+{\gamma_1}-\a_1}^{g+h+k}\neq0$ and ${\b_1}+{\gamma_1}-\a_1\in\Lam^{g+h+k}.$ We get that $\{\a_1, -\gamma_1, {\b_1}+{\gamma_1}-\a_1\}$ is a connection from $\a_1$ to $\b_1$ and so $\a_1\sim\b_1,$ which is a contradiction.. Hence,
\begin{equation}\label{50'.}
\rho(\LL_{\b_1}^h, \LL_{\gamma_1}^k)(A_{-\a_1}^g)\LL_{\a_1}^{g'}=0.
\end{equation}
It follows 
\begin{equation}\label{51}
[I_{0, [\a]}, \vv_{[\b]}, \vv_{[\gamma]}]=0.
\end{equation}
By a similar argument, we also get, 
\begin{equation}\label{51'}
[\vv_{[\a]}, I_{0, [\b]}, \vv_{[\gamma]}]=[\vv_{[\a]}, \vv_{[\b]}, I_{0, [\gamma]}]=0.
\end{equation}

 Finally, we consider the summand $[I_{0, [\a]}, I_{0, [\b]}, \vv_{[\gamma]}]$ in Eqs. (\ref{205}). One can write,
 \begin{eqnarray*}
 [I_{0, [\a]}, I_{0, [\b]}, \vv_{[\gamma]}] &=&[(\sum_{\substack{{\a_1}\in[\a]\\ {-\a_1}\in\Lam}} A_{-\a_1}^g\LL_{\a_1}^{g'})+(\sum_{\substack{{\a_1}, {\a_2}, {\a_3}\in[\a]\\ {\a_1}+{\a_2}+{\a_3}=0}}[\LL_{\a_1}^{g_1},
 \LL_{\a_2}^{g_2}, \LL_{\a_3}^{g_3}]),\\
 &~& (\sum_{\substack{{\b_1}\in[\b]\\ {-\b_1}\in\Lam}} A_{-\b_1}^h\LL_{\b_1}^{h'})+(\sum_{\substack{{\b_1}, {\b_2}, {\b_3}\in[\b]\\ {\b_1}+{\b_2}+{\b_3}=0}}[\LL_{\b_1}^{h_1}, \LL_{\b_2}^{h_2}, \LL_{\b_3}^{h_3}]), \vv_{[\gamma]}].
 \end{eqnarray*}
 The above statement includs four items which we consider in the following. First, consider the item $[A_{-\a_1}^g\LL_{\a_1}^{g'},  A_{-\b_1}^h\LL_{\b_1}^{h'}, \LL_{\gamma_1}^k],$ where ${\gamma_1}\in[\gamma].$ By  Eqs. (\ref{4}), (\ref{149}) and (\ref{50'.}), we have
 
\begin{eqnarray*}
 	[A_{-\a_1}^g\LL_{\a_1}^{g'},  A_{-\b_1}^h\LL_{\b_1}^{h'}, \LL_{\gamma_1}^k] &=&[\LL_{\gamma_1}^k, A_{-\a_1}^g\LL_{\a_1}^{g'},  A_{-\b_1}^h\LL_{\b_1}^{h'}]\\ 
 	&=& A_{-\b_1}^h[\LL_{\gamma_1}^k, A_{-\a_1}^g\LL_{\a_1}^{g_1}, \LL_{\b_1}^{h_1}]\\
 	&+&A_{-\a_1}^{g}\rho(\LL_{\gamma_1}^k, \LL_{\a_1}^{g_1})(A_{-\b_1}^h)\LL_{\a_1}^{g_1}\\
 	&=&A_{-\b_1}^hA_{-\a_1}g[\LL_{\gamma_1}^k, \LL_{\a_1}^{g_1}, \LL_{\b_1}^{h_1}]\\
 	&+&A_{-\b_1}^h\rho(\LL_{\b_1}^{h_1}, \LL_{\gamma_1}^k)(A_{-\a_1}^g)\LL_{\a_1}^{g_1}\\
 	&+&A_{-\a_1}^g\rho(\LL_{\gamma_1}^k, \LL_{\a_1}^{g_1})(A_{-\b_1}^h)\LL_{\b_1}^{h_1}\\
 	&=&0.
\end{eqnarray*}
Next, by  Eq. (\ref{149}) we have
\begin{equation*}
[A_{-\a_1}^g\LL_{\a_1}^{g_1}, [\LL_{\b_1}^{h_1}, \LL_{\b_2}^{h_2}, \LL_{\b_3}^{h_3}], \LL_{\gamma_1}^k]\subset[\vv_{[\a]}, \vv_{[\b]}, \vv_{[\gamma]}]=0,	
\end{equation*}
and
\begin{equation*}
[[\LL_{\a_1}^{g_1},
\LL_{\a_2}^{g_2}, \LL_{\a_3}^{g_3}],  A_{-\b_1}^h\LL_{\b_1}^{h'}, \LL_{\gamma_1}^k]\subset[\vv_{[\a]}, \vv_{[\b]}, \vv_{[\gamma]}]=0,	
\end{equation*}
and
\begin{equation*}
[[\LL_{\a_1}^{g_1}, \LL_{\a_2}^{g_2}, \LL_{\a_3}^{g_3}],  [\LL_{\b_1}^{h_1}, \LL_{\b_2}^{h_2}, \LL_{\b_3}^{h_3}], \LL_{\gamma_1}^k]\subset[\vv_{[\a]}, \vv_{[\b]}, \vv_{[\gamma]}]=0.	
\end{equation*}
Summarizing a discussion of above, we get
$$
	[I_{0, [\a]}, I_{0, [\b]}, \vv_{[\gamma]}]=0
$$
In a similar way we get
$$
[I_{0, [\a]}, \vv_{[\b]}, I_{0,
	[\gamma]}]=0,~~	[\vv_{[\a]}, I_{0, [\b]}, \vv_{[\gamma]}]=0.
$$
 Therfore,
 $$
[I_{[\a]}, I_{[\b]}, I_{[\gamma]}]=0
$$
By a similar argument as above one can prove $[I_{[\a]}, I_{[\a]}, I_{[\b]}]=0.$\qed

\begin{thm}\label{main1} The following assertions hold
\begin{itemize}
\item[(1)] For any $[\a]\in\Pi/\sim,$ the linear graded subspace
$$
I_{[\a]} =I_{0, [\a]}\oplus\vv_{[\a]},
$$
of split $3-$Lie-Rinehart color algebra $(\LL, A)$ associated to $[\a]$ is an  ideal of $(\LL, A).$
\item[(2)]  If $(\LL, A)$ is simple, then there exists a connection from
$\a$ to $\b$ for any $\a, \b\in\Pi$ and
$$
\hh=(\sum_{\substack{\b\in\Pi, -\b\in\Lam\\ g, g'\in G}} A_{-\b}^g\LL_\b^{g'})
+(\sum_{\substack{\b, \gamma, \mu\in\Pi\\ \b+\gamma+\mu=0\\ g, g', g''\in G}}[\LL_\b^g,
\LL_{\gamma}^{g'}, \LL_\mu^{g''}]).
$$
\end{itemize}
\end{thm}
\noindent {\bf Proof.}

(1) Since
$$
[I_{[\a]}, \hh, \hh]=[I_{[\a]}, \LL_0, \LL_0]\subset\vv_{[\a]},
$$ 
By Propositions \ref{subalg} and \ref{12345},  we have
\begin{eqnarray*}\label{2080}
[I_{[\a]}, \LL, \LL]&=&[I_{[\a]},
\hh\oplus(\bigoplus_{\b\in[\a]}\LL_\b)\oplus(\bigoplus_{\gamma\notin[\a]}\LL_\gamma), \hh\oplus(\bigoplus_{\b\in[\a]}\LL_\b)\oplus(\bigoplus_{\gamma\notin[\a]}\LL_\gamma)]\\
&\subset&[I_{[\a]}, \hh, \LL_\b]+[I_{[\a]}, \hh, \LL_\gamma]+[I_{[\a]},, \LL_\a, \LL_\a]+[I_{[\a]}, \LL_\a, \LL_\gamma]\\
&+&[I_{[\a]}, \LL_\gamma, \hh]+[I_{[\a]}, \LL_\gamma, \LL_\a]+[I_{[\a]}, \LL_\gamma, \LL_\gamma]\\
&\subset&I_{[\a]}.
\end{eqnarray*}
Then $I_{[\a]}$ is a $3-$Lie ideal of $\LL.$ We also have $A I_{[\a]}\subset I_{[\a]},$ thanks to Proposition \ref{subalg}-(2), and
$$
\rho( I_{[\a]}, \LL)(A)\LL\subset A[I_{[\a]}, \LL, \LL]+[I_{[\a]}, \LL, A\LL]I_{[\a]}\subset I_{[\a]}.
$$
Follows from Proposition \ref{subalg}, $I_{[\a]}$is an ideal of $(\LL, A).$

(2)  The simplicity of $(\LL, A)$ implies that $I_{[\a]}\in\{\LL, \ker\rho\}.$ If $I_{[\a]}=\LL$ for some $\a\in\Pi,$ we are done. Otherwise, if $I_{[\a]}=\ker\rho$ for all $\a\in\Pi$ we have $[\a]=[\b]$ for any $\b\in\Pi$ and again $\Pi=[\a].$ We  conclude that $\LL$ has all of its non zero roots connected and therefore, by Eqs. (\ref{119.}), (\ref{38'}) and (\ref{38''}), we have
$$
\hh=(\sum_{\substack{\b\in\Pi, -\b\in\Lam\\ g, g'\in G}} A_{-\b}^g\LL_\b^{g'})
+(\sum_{\substack{\b, \gamma, \mu\in\Pi\\ \b+\gamma+\mu=0\\ g, g', g''\in G}}[\LL_\b^g,
\LL_{\gamma}^{g'}, \LL_\mu^{g''}]).
$$
\qed
\begin{thm}\label{main2} For a vector space complement $\uu$ in $\hh$ of
$$
(\sum_{\substack{\b\in\Pi, -\b\in\Lam\\ g, g'\in G}} A_{-\b}^g\LL_\b^{g'})
+(\sum_{\substack{\b, \gamma, \mu\in\Pi\\ \b+\gamma+\mu=0\\ g, g', g''\in G}}[\LL_\b^g,
\LL_{\gamma}^{g'}, \LL_\mu^{g''}]),
$$  
we have
$$
\LL=\uu\oplus\sum_{[\a]\in\Pi/\sim}I_{[\a]},
$$
where any $I_{[\a]}$ is one of the graded ideals of $(\LL, A)$
described in Theorem \ref{main1}-(1), and satisfying the conditions in  Proposition \ref{12345}. 
\end{thm}
\noindent {\bf Proof.} Each $I_{[\a]}$ is well defined and,  by
Theorem \ref{main1}-(1), a graded ideal of $(\LL, A).$ It is
clear that
$$
\LL=\hh\oplus(\bigoplus_{\a\in\Pi}\LL_\a)=\uu\oplus\sum_{[\a]\in\Pi/\sim}I_{[\a]}.
$$
Finally, Proposition \ref{12345} gives us
$$
[I_{[\a]}, I_{[\b]}, I_{[\gamma]}]=0,~~~\hbox{and~~}~~~~[I_{[\a]}, I_{[\a]}, I_{[\b]}]=0,
$$
whenever $[\a], [\b], [\gamma]\in\Pi/\sim$ be different from each other.\qed

\begin{pro}\label{2.17} If $Z_\rho(\LL)=\{0\}$ and 
$$
\hh=(\sum_{\substack{\b\in\Pi, -\b\in\Lam\\ g, g'\in G}} A_{-\b}^g\LL_\b^{g'})
+(\sum_{\substack{\b, \gamma, \mu\in\Pi\\ \b+\gamma+\mu=0\\ g, g', g''\in G}}[\LL_\b^g,
\LL_{\gamma}^{g'}, \LL_\mu^{g''}]),
$$	
then $\LL$ is the direct sum of the graded ideals
given in Theorem \ref{main1}-(1),
$$
\LL=\bigoplus_{[\a]\in\Pi/\sim}I_{[\a]},
$$	
satisfying the conditions in  Proposition \ref{12345}.
\end{pro}
\noindent {\bf Proof.} From Theorem \ref{main2} and
$$
\hh=(\sum_{\substack{\b\in\Pi, -\b\in\Lam\\ g, g'\in G}} A_{-\b}^g\LL_\b^{g'})
+(\sum_{\substack{\b, \gamma, \mu\in\Pi\\ \b+\gamma+\mu=0\\ g, g', g''\in G}}[\LL_\b^g,
\LL_{\gamma}^{g'}, \LL_\mu^{g''}]),
$$
it is clear that
$\LL=\sum_{[\a]\in\Pi/\sim}I_{[\a]}.$ For the direct character, take some $x\in I_{[\a]}\cap\sum_{[\a]\in\Pi/\sim, [\b]\neq [\a]}I_{[\b]}.$ By Proposiyion \ref{12345}, we have 

\begin{eqnarray*}
[x, \LL, \LL]&\subset&[x, I_{[\a]}, I_{[\a]}]+[x, I_{[\a]}, \sum_{\b\notin [\a]}I_{[\b]}]\\
&+&[x, \sum_{\b\in [\a]}I_{[\b]}, \sum_{\b\notin [\a]}I_{[\b]}]=0.
\end{eqnarray*}
That is, $x\in Ann(\LL).$ We also have $\rho(x, \LL)=0,$ and so $x\in Z_\rho(\LL)=\{0\}.$
Hence, $\LL=\bigoplus_{[\a]\in\Pi/\sim}I_{[\a]}$,  as by Proposition \ref{12345},  
$$
[I_{[\a]}, I_{[\b]}, I_{[\gamma]}]=0,~~~\hbox{and~~}~~~~[I_{[\a]}, I_{[\a]}, I_{[\b]}]=0,
$$
whenever $[\a], [\b], [\gamma]\in\Pi/\sim$ be different from each other.\qed

\section{ CONNECTIONS IN THE WEIGHT SYSTEM OF A. DECOMPOSITION of A}
\setcounter{equation}{0}\

We begin this section by introducing an adequate notion of connection among the
weights in $\Lam.$  For a split regular  $3-$Lie-Rinehart color algebra $(\LL, A)$, since $A$ is a commutative
associative color algebra, then the decomposition of $A$ is similar to \cite{ABCS} and omit the proof of some results.

\begin{DEF}\label{conn}
	Let $\lam, ~\mu\in \Lam,$ we say that {\em $\lam$
		is connected to $\mu$} and denoted by $\lam\approx\mu,$  if either $\lam=\pm\mu,$ or there exists a
	family
	$$
	\{\lam_1, \lam_2, \lam_3, ..., \lam_k\}\subset\pm\Pi\cup\pm\Lam,
	$$
with $k\geq 2,$ such that	satisfying the following conditions;\\
\begin{itemize}	
\item[(1)] $\lam_1=\lam.$
	
\item[(2)] $\lam_1+\lam_2\in\pm\Lam,\\
		\lam_1+ \lam_2+ \lam_3\in\pm\Lam,\\
		...\\
		\lam_1+ \lam_2+ \lam_3+...+\lam_{k-1}\in\pm\Lam.$
		
\item[(3)] $	\lam_1+ \lam_2+ \lam_3+...+\lam_k\in\{\mu, -\mu\}.$
	\end{itemize}
	The family $\{\lam_1, \lam_2, \lam_3, ..., \lam_k\}$ is called a
	{\em connection} from $\lam$ to $\mu.$
\end{DEF}
 
As in the previous section we can prove the next result;
\begin{pro}\label{rela2}
	The relation $\approx$ in $\Lam$ defined by
	$$
	\lam\approx\mu~\hbox{ if and only if}~ \lam~\hbox{ is connected to}~ \mu,
	$$
	is an equivalence relation.
\end{pro}

\begin{rem}\label{lem}
Let $\lam, \mu\in\Lam$ such that $\lam\approx\mu.$ If $\lam+\eta\in\Lam,$ for $\eta\in\Pi\cup\Lam$ then $\lam\approx\mu+\eta.$ Considering the connection $\{\mu, \eta\}$ we get $\mu\approx\mu+\eta$ and by transitivity $\lam\approx\mu+\eta.$
\end{rem}

By the  Proposition \ref{rela2}, we can consider the equivalence relation
in $\Lam$ by the connection relation $\approx$ in $\Lam.$ So we denote
by
$$
\Lam/\approx :=\{[\lam] : \lam\in\Lam\},
$$
where $[\lam]$ denotes  the set of non zero weights  which are
connected to $\lam.$ 

Our next goal in this section is to associate an adequate ideal
$\aa_{[\lam]}$ of $A$ to any $[\lam].$ For a fixed  $\lam\in\Lam,$ we define
$$
A_{0, [\lam]} :=(\sum_{\mu\in[\lam]} A_{-\mu}A_\mu)+(\sum_{\substack{\a, \b\in\Pi, \mu\in[\lam]\\ \a+\b+\mu=0}}\rho(\LL_\a,
\LL_\b) A_\mu)\subset A_0.
$$
That is
\begin{eqnarray}\label{119}
	A_{0, [\lam]} : =(\sum_{\substack{\mu\in [\lam]\\  g\in G}} A_{-\mu}^{-g}A_\mu^{g})
	+(\sum_{\substack{\a, \b\in\Pi, \mu\in[\lam]\\ \a+\b+\mu=0\\ g, g', g''\in G}}\rho(\LL_\a^g,
	\LL_\b^{g'}) A_\mu^{g''})\subset A_0^{g+g'+g''}.
\end{eqnarray}
Next, we define
\begin{equation}\label{48'}
	A_{[\lam]}
	:=\bigoplus_{\mu\in[\lam]}A_\mu=\bigoplus_{g\in G}\bigoplus_{\mu\in[\lam]}A^g_\mu.
\end{equation}
Finally, we denote by $\aa_{[\lam]}$ the direct sum of the two graded
subspaces above, that is,
\begin{equation}\label{48''}
	\aa_{[\lam]} :=A_{0, [\lam]}\oplus A_{[\lam]}.
\end{equation}

\begin{pro}\label{4.7} For any $\lam, \mu\in\Lam,$ the following assertions hold.
\begin{itemize}
\item[(1)] $\aa_{[\lam]} \aa_{[\lam]}\subset\aa_{[\lam]}.$

\item[(2)]  If $[\lam]\neq [\mu],$ then $\aa_{[\lam]} \aa_{[\mu]}=0.$
\end{itemize} 
\end{pro}
\noindent {\bf Proof.} We only prove (1) and similar for (2). By Eq. (\ref{48''}) and comutativity of $A$ we have

\begin{equation}\label{49}
\aa_{[\lam]} \aa_{[\lam]}=(A_{0, [\lam]}\oplus A_{[\lam]})(A_{0, [\lam]}\oplus A_{[\lam]})\subset	A_{0, [\lam]} A_{0, [\lam]}+A_{0, [\lam]}A_{[\lam]}+A_{[\lam]} A_{[\lam]}.
\end{equation}
Consider the second summand in Eq. (\ref{49}). Given $\mu\in [\lam], g\in G$ we have 
$$
A_{0, [\lam]}A_\mu^g\subset A_0 A^g_\mu\subset A^g_\mu.
$$
Therefore
\begin{equation}\label{50.}
A_{0, [\lam]}A_{[\lam]}\subset \aa_{[\lam]}.
\end{equation}
Let us consider the third summand in Eq. (\ref{49}). Given $\mu, \eta\in[\Lam], g, g'\in G$ such that $A_\mu^g A_\eta^{g'}\neq 0.$ Then $A_\mu^g A_\eta^{g'}\subset A_{\mu+\eta}^{g+g'}.$ If $\mu+\eta=0,$ by Eq. (\ref{119}) we have $A_\mu^g A_{-\mu}^{g'}\subset A_{0, [\lam]}.$ Suppose $\mu+\eta\in\Lam^{g+g'},$ then thanks to Remark \ref{lem}, we have $\mu+\eta\in[\lam]$ and so $A_\mu^g A_{\eta}^{g'}\subset A_{\mu+\eta}^{g+g'}\subset A_{[\lam]}.$ Hence,
$$
(\bigoplus_{g\in G}\bigoplus_{\mu\in[\lam]}A^g_\mu)(\bigoplus_{{g'}\in G}\bigoplus_{\eta\in[\lam]}A^{g'}_\eta)\subset A_{0, [\lam]}\oplus A_{[\lam]},
$$
and so
\begin{equation}\label{51.}
A_{[\lam]} A_{[\lam]}\subset \aa_{[\lam]}.
\end{equation}
Finally, consider the first summand in Eq. (\ref{49}). Given $\mu_i, \eta_i\in\Lam, \a_i, \b_i\in\Pi, i=1, 2$ and $g, g_1, g'_1, g''_1, h, h_2, h'_2, h''_2\in G$ such that

\begin{eqnarray}\label{52}
\nonumber 0&\neq&	 (A_{-\mu_1}^{-g} A_{\mu_1}^g+ \rho(\LL_{\a_1}^{g_1},
\LL_{\b_1}^{g'_1})( A_{\eta_1}^{g''_1}))(A_{-\mu_2}^{-h} A_{\mu_2}^{h} + \rho(\LL_{\a_2}^{h_2},
\LL_{\b_2}^{h'_2}) (A_{\eta_2}^{h''_2}))\\
&=& (A_{-\mu_1}^{-g} A_{\mu_1}^g)(A_{-\mu_2}^{-h} A_{\mu_2}^{h})+\rho(\LL_{\a_1}^{g_1},
\LL_{\b_1}^{g'_1}) (A_{\eta_1}^{g''_1})(A_{-\mu_2}^{-h} A_{\mu_2}^{h})\\
\nonumber&+& (A_{-\mu_1}^{-g} A_{\mu_1}^g) \rho(\LL_{\a_2}^{h_2},
\LL_{\b_2}^{h'_2})( A_{\eta_2}^{h''_2})+ \rho(\LL_{\a_1}^{g_1},
\LL_{\b_1}^{g'_1}) (A_{\eta_1}^{g''_1})\rho(\LL_{\a_2}^{h_2},
\LL_{\b_2}^{h'_2})( A_{\eta_2}^{h''_2}).
\end{eqnarray}
We are going to consider all sammands in Eq. (\ref{52}), in four cases; 

{\bf Case 1.} For the first sammand, if $\mu_1+\mu_2\neq 0,$ by Remark \ref{lem} and the fact that $A$ is a commutative and associative algebra, we have  
\begin{equation}\label{53.}
 ~~~~(A_{-\mu_1}^{-g} A_{\mu_1}^g)(A_{-\mu_2}^{-h} A_{\mu_2}^{h})=(A_{-\mu_1}^{-g}A_{-\mu_2}^{-h})( A_{\mu_1}^gA_{\mu_2}^{h})\subset A_{-(\mu_1+\mu_2)}^{-(g+h)} A_{(\mu_1+\mu_2)}^{g+h}\subset  A_{0, [\lam]}.
\end{equation}
If $\mu_1+\mu_2= 0,$ it follows
\begin{equation}\label{54.}
	~~~~(A_{-\mu_1}^{-g} A_{\mu_1}^g)(A_{-\mu_2}^{-h} A_{\mu_2}^{h})=A_{-\mu_1}^{-g}(A_{-\mu_2}^{-h} A_{\mu_2}^h)A_{\mu_1}^{g}\subset A_{-\mu_1}^{-g} A_{\mu_1}^{g}\subset  A_{0, [\lam]}.
\end{equation}
From Eqs. (\ref{53.}) and (\ref{54.}), we get 

\begin{equation}\label{55}
(A_{-\mu_1}^{-g} A_{\mu_1}^g)(A_{-\mu_2}^{-h} A_{\mu_2}^{h})\subset \aa_{[\lam]}.
\end{equation}

{\bf Case 2.}  For the second sammand, since $\rho(\LL_{\a_1}^{g_1},
\LL_{\b_1}^{g'_1})$ is a derivation in $A$ we have  
\begin{eqnarray*}
\rho(\LL_{\a_1}^{g_1}, \LL_{\b_1}^{g'_1}) (A_{\eta_1}^{g''_1})(A_{-\mu_2}^{-h} A_{\mu_2}^{h})&\subset& \rho(\LL_{\a_1}^{g_1}, \LL_{\b_1}^{g'_1}) (A_{\eta_1}^{g''_1}(A_{-\mu_2}^{-h} A_{\mu_2}^{h}))\\
&+& A_{\eta_1}^{g''_1}\rho(\LL_{\a_1}^{g_1}, \LL_{\b_1}^{g'_1})(A_{-\mu_2}^{-h} A_{\mu_2}^{h}).
\end{eqnarray*}
Thanks to  Proposition \ref{2.7}, 
$$
\rho(\LL_{\a_1}^{g_1}, \LL_{\b_1}^{g'_1}) (A_{\eta_1}^{g''_1}(A_{-\mu_2}^{-h} A_{\mu_2}^{h}))\subset \rho(\LL_{\a_1}^{g_1}, \LL_{\b_1}^{g'_1}) (A_{\eta_1}^{g''_1})\subset A_{\a_1+\b_1+\eta_1}^{g_1+g'_1+g''_1},
$$
and
$$
A_{\eta_1}^{g''_1}\rho(\LL_{\a_1}^{g_1}, \LL_{\b_1}^{g'_1})(A_{-\mu_2}^{-h} A_{\mu_2}^{h})\subset A_{\eta_1}^{g''_1} A_{\a_1+\b_1}^{g_1+g'_1}\subset A_{\a_1+\b_1+\eta_1}^{g_1+g'_1+g''_1}.
$$
Therefore, we get 
\begin{equation}\label{56.}
\rho(\LL_{\a_1}^{g_1}, \LL_{\b_1}^{g'_1}) (A_{\eta_1}^{g''_1})(A_{-\mu_2}^{-h} A_{\mu_2}^{h})\subset \aa_{[\lam]}.
\end{equation}
{\bf Case 3.} Consider the third summand in Eq. (\ref{52}),by the commtativity of $A,$ we have
\begin{equation*}
(A_{-\mu_1}^{-g} A_{\mu_1}^g) \rho(\LL_{\a_2}^{h_2},
\LL_{\b_2}^{h'_2})( A_{\eta_2}^{h''_2})\subset A_{\a_2+\b_2+\eta_2}^{h_2+h'_2+h''_2}\subset\aa_{[\lam]}. 
\end{equation*}
{\bf Case 4.} For the last summand in Eq. (\ref{52}), since $\rho(\LL_{\a_1}^{g_1},
\LL_{\b_1}^{g'_1})$ is a derivation in $A,$ we have  
\begin{eqnarray*}
 \rho(\LL_{\a_1}^{g_1},
\LL_{\b_1}^{g'_1}) (A_{\eta_1}^{g''_1})\rho(\LL_{\a_2}^{h_2},
\LL_{\b_2}^{h'_2})( A_{\eta_2}^{h''_2})	&\subset&  \rho(\LL_{\a_1}^{g_1},
\LL_{\b_1}^{g'_1}) (A_{\eta_1}^{g''_1}(\rho(\LL_{\a_2}^{h_2},
\LL_{\b_2}^{h'_2})( A_{\eta_2}^{h''_2}))\\
&+&A_{\eta_1}^{g''_1} \rho(\LL_{\a_1}^{g_1},
\LL_{\b_1}^{g'_1})(\rho(\LL_{\a_2}^{h_2},
\LL_{\b_2}^{h'_2})( A_{\eta_2}^{h''_2})).
\end{eqnarray*}
Now by Proposition \ref{2.7}, we have
\begin{eqnarray*}
\rho(\LL_{\a_1}^{g_1},
\LL_{\b_1}^{g'_1}) (A_{\eta_1}^{g''_1}(\rho(\LL_{\a_2}^{h_2},
\LL_{\b_2}^{h'_2})( A_{\eta_2}^{h''_2}))&\subset&\rho(\LL_{\a_1}^{g_1},
\LL_{\b_1}^{g'_1})(A_{\a_2+\b_2+\eta_1+\eta_2}^{g''_1+h_2+h'_2+h''_2})\\
&\subset& A_{\a_1+\b_1+\a_2+\b_2+\eta_1+\eta_2}^{g_1+g'_1+g''_1+h_2+h'_2+h''_2}.
\end{eqnarray*}
Similarly,
$$
A_{\eta_1}^{g''_1} \rho(\LL_{\a_1}^{g_1},
\LL_{\b_1}^{g'_1})(\rho(\LL_{\a_2}^{h_2},
\LL_{\b_2}^{h'_2})( A_{\eta_2}^{h''_2}))\subset A_{\a_1+\b_1+\a_2+\b_2+\eta_1+\eta_2}^{g_1+g'_1+g''_1+h_2+h'_2+h''_2}.
$$
Therefore,
\begin{equation*}
 \rho(\LL_{\a_1}^{g_1},
\LL_{\b_1}^{g'_1}) (A_{\eta_1}^{g''_1})\rho(\LL_{\a_2}^{h_2},
\LL_{\b_2}^{h'_2})( A_{\eta_2}^{h''_2})\subset\aa_{[\lam]}. 
\end{equation*}
Summarizing a discussion of above four cases, we get
\begin{equation}\label{58.}
A_{0, [\lam]} A_{0, [\lam]}\subset\aa_{[\lam]}. 
\end{equation}
From Eqs. (\ref{50.}), (\ref{51.}) and (\ref{58.}) we conclude the result.\qed

We recall that a $G-$graded subspace $I$ of a commutative and associative color algebra $A$ is called an ideal of $A$ if $A I \subset I.$ We say that $A$ is simple if $A A\neq 0$ and it contains no proper ideals.

\begin{thm}\label{main1'} Let $A$ be a commutative and associative color algebra associated to a  split $3-$Lie-Rinehart color algebra $(\LL, A).$ Then the following assertions hold.
	\begin{itemize}
		\item[(1)] For any $[\lam]\in\Lam/\approx,$ the linear graded subspace
		$$
		\aa_{[\lam]} =A_{0, [\lam]}\oplus A_{[\lam]},
		$$
		of  color algebra $A$ associated to $[\lam]$ is an  ideal of $A.$
		\item[(2)]  If $A$  is simple then all weights of $\Lam$ are connected. Furthermore,
		$$
		A_0=(\sum_{\substack{\mu\in\Lam\\  g\in G}} A_\mu^{g} A_{-\mu}^{-g})
		+(\sum_{\substack{\a, \b\in\Pi, \mu\in\Lam\\ \a+\b+\mu=0\\ g, g', g''\in G}}\rho(\LL_\a^g,
		\LL_\b^{g'}) A_\mu^{g''}).
		$$
	\end{itemize}
\end{thm}
\noindent {\bf Proof.} (1) By Eq. (\ref{00.}), we have
$$
\aa_{[\lam]} A=\aa_{[\lam]}(A^0\oplus(\bigoplus_{g\in G}(\bigoplus_{\mu\in[\Lam]^g}A^g_\mu+\bigoplus_{\mu\notin[\Lam]^g}A^g_\mu))).
$$
Now, by asociativity of $A$, we have $\aa_{[\lam]}A^0\subset \aa_{[\lam]},$ and by Proposition \ref{4.7}, 
$$
\aa_{[\lam]}(\bigoplus_{g\in G}(\bigoplus_{\mu\in[\Lam]^g}A^g_\mu))\subset \aa_{[\lam]},
$$ 
and 
$$
\aa_{[\lam]}(\bigoplus_{g\in G}(\bigoplus_{\mu\notin[\Lam]^g}A^g_\mu))\subset \aa_{[\lam]}.
$$ 
Hence,
$$
\aa_{[\lam]} A\subset \aa_{[\lam]}.			
$$
That is, $\aa_{[\lam]}$ is an  ideal of $A.$

(2) The simplicity of $A$ implies $\aa_{[\lam]}=A,$ for any $\lam\in\Lam.$ From here it is clear that $\Lam=[\lam],$ and so $	A_0=(\sum_{\substack{\mu\in\Lam\\  g\in G}} A_\mu^{g} A_{-\mu}^{-g})
+(\sum_{\substack{\a, \b\in\Pi, \mu\in\Lam\\ \a+\b+\mu=0\\ g, g', g''\in G}}\rho(\LL_\a^g,
\LL_\b^{g'}) A_\mu^{g''}).$\qed

\begin{thm}\label{main2'} Let $A$ be a commutative and associative color algebra associated to a  split $3-$Lie-Rinehart color algebra $(\LL, A).$ Then 
$$
A=\v+\sum_{[\lam]\in\Lam/\approx}\aa_{[\lam]},
$$
where $\v$ is a graded linear complement in $A^0$ of 
$$
(\sum_{\substack{\mu\in\Lam\\  g\in G}} A_\mu^{g} A_{-\mu}^{-g})
+(\sum_{\substack{\a, \b\in\Pi, \mu\in\Lam\\ \a+\b+\mu=0\\ g, g', g''\in G}}\rho(\LL_\a^g,
\LL_\b^{g'}) A_\mu^{g''}),
$$ 
and any $\aa_{[\lam]}$ is one of ideals of $A$ described in Theorem \ref{main1'}-(1). Forthermore  $\aa_{[\lam]} \aa_{[\mu]}=0,$ when $[\lam]\neq [\mu].$
\end{thm}
\noindent {\bf Proof.} By Theorem \ref{main1'}-(1), $\aa_{[\lam]}$ is a well defined ideal of $A,$ bing clear that 
$$
A=\bigoplus_{g\in G}(A^g\oplus(\bigoplus_{\lam\in\Lam}A^g_\lam))=A^0\oplus(\bigoplus_{g\in G}\bigoplus_{\lam\in\Lam^g}A^g_\lam)=\v+\sum_{[\lam]\in\Lam/\approx}\aa_{[\lam]},
$$
and by Proposition \ref{4.7}-(2), if  $[\lam]\neq [\mu]$ then $\aa_{[\lam]} \aa_{[\mu]}=0.$\qed

Let us denote by $Ann(A) := \{a \in A : a A = 0\}$ the annihilator of the commutative and
associative algebra $A.$

\begin{cor}\label{4.19} Let $(\LL, A)$ be a $3-$Lie-Rinehart colr algebra. If $Ann(A) = 0$ and
$$
A_0=(\sum_{\substack{\mu\in\Lam\\  g\in G}} A_\mu^{g} A_{-\mu}^{-g})
+(\sum_{\substack{\a, \b\in\Pi, \mu\in\Lam\\ \a+\b+\mu=0\\ g, g', g''\in G}}\rho(\LL_\a^g,
\LL_\b^{g'}) A_\mu^{g''}),
$$
then $A$ is the direct sum of the ideals given in Theorem \ref{main1'}-(1),
$$
A=\bigoplus_{[\lam]\in\Lam/\approx}\aa_{[\lam]}.
$$
Furthermore, $\aa_{[\lam]} \aa_{[\mu]}=0,$ when $[\lam]\neq [\mu].$
\end{cor}
\noindent {\bf Proof.} This can be proved analogously to Corollary 3.8 in \cite{ABCS}.

\section{Relating the dcompositions of $\LL$ and $A$} 
\setcounter{equation}{0}\

In this section,  we will show that the decompositions of $\LL$ and $A$ as direct sum of
ideals, given in Sections 2 and 3 respectively, are closely related.

\begin{lem}  Let $(\LL, A)$ be a  split $3-$Lie-Rinehart color algebra and $I$ an ideal of $\LL$.
Then $I = (I\cap\hh) \oplus(\bigoplus_{\a\in\Pi}(I\cap\LL_\a))$
\end{lem}
\noindent {\bf Proof.}  Since $(\LL, A)$ is split, we get $\LL = \hh\oplus (\bigoplus_{\a\in\Pi}\LL_\a).$ By the assumption that $I$ is an
ideal of L, it is clear that $I$ is a submodule of $\LL$. Thus $I$ is a weight module and therefore
$I = (I\cap\hh) \oplus(\bigoplus_{\a\in\Pi}(I\cap\LL_\a)).$\qed

Observe that for  a  split $3-$Lie-Rinehart color algebra, by grading, let us assert
that given any non-zero graded ideal $I$ of $\LL$ we can write
\begin{equation}\label{301}
I=\bigoplus_{g\in G}((I\cap\hh^g)\oplus(\bigoplus_{\a\in\Pi^g_I}(I\cap\LL_\a^g))),
\end{equation}
where $\Pi^g_I :=\{\a\in\Pi : I\cap\LL_\a^g\neq0\}$
for each $g\in G.$

\begin{lem}\label{5.3}
 Let $(\LL, A)$ be a  split $3-$Lie-Rinehart color algebra  with $Z_\rho(\LL) = 0$ and
$I$ an ideal of $\LL$. If $I \subset\hh,$ then $I = \{0\}.$
\end{lem}
\noindent {\bf Proof.} By Proposition \ref{2.17} and $I \subset\hh,$ we have $I\cap\LL_a=0,$ for any $\a\in\Pi.$ So
$$
[I, \hh, \LL_\a]=0,~~[I, \LL_\a, \LL_\b]=0, \forall \a, \b\in\Pi.
$$
Therefore,
$$
[I, \LL, \LL] = [I, \hh\oplus\bigoplus_{\a\in\Pi}\LL_\a, \hh\oplus\bigoplus_{\b\in\Pi}\LL_\b]\subset [I, \hh, \bigoplus_{\b\in\Pi}\LL_\b]+[I, \bigoplus_{\a\in\Pi}\LL_\a, \bigoplus_{\b\in\Pi}\LL_\b]=0.
$$
So $I \subset Z_\rho(\LL) = 0.$\qed

\begin{pro}\label{5.4}
If a split $3-$Lie-Rinehart color algebra $(\LL, A)$ can be decomposed
into a direct sum of finite ideals $\LL_i
, 1 \leq i \leq s,$ and $Z_{\LL}(A) = 0,$ then  $(\LL, A)$ is a
split $3-$Lie-Rinehart color algebra with a splitting Cartan subalgebra $\hh$  if and only
if $(\LL_i, A, \rho|_{\LL_i\times\LL_i}), 1 \leq i \leq s$ are split $3-$Lie-Rinehart color algebras with a splitting Cartan subalgebra $\hh_i,$ respectively, such that $\hh =\bigoplus_{i=1}^s\hh_i,~~\Pi=\bigcup_{i=1}^s\Pi_i,$ and $\Lam=\bigcup_{i=1}^s\Lam_i,$ where $\Pi_i$ are
root systems of $\LL_i$ and $\Lam_i$ are weight system of $A$ associated to $\hh_i,$ respectively.
\end{pro}
\noindent {\bf Proof.} Suppose $(\LL, A)$ is a split $3-$Lie-Rinehart color algebra with Cartan subalgebra $\hh$ and root system $\Pi.$ Let $\LL=\bigoplus_{i=1}^s\LL_i,$ where $\LL_i$ be non-zero ideals of $\LL.$ By Eq. ( \ref{301}),
\begin{equation}
\LL_i=\bigoplus_{g\in G}((\LL_i\cap\hh^g)\oplus(\bigoplus_{\a\in\Pi}(\LL_i\cap\LL_\a^g))),~i=1, 2, 3, ..., s.
\end{equation}
If $ \LL_i\cap\LL_\a^g\neq0,~~\forall g\in G$ then $\a|_{\hh_i\times\hh_i}\neq 0,$ where $\hh_i=\bigoplus_{g\in G}(\LL_i\cap\hh^g)=\bigoplus_{g\in G}\hh_i^g.$ In fact, if  $\a|_{\hh_i\times\hh_i}= 0.$ For any $0\neq x_i\in\LL_i\cap\LL_\a^g,$ from $\hh=\bigoplus_{i=1}^s\hh_i=\bigoplus_{i=1}^s(\bigoplus_{g\in G}\hh_i^g),$ we have
$$
[\hh, \hh, x_i]\subset \sum_{i=1}^s(\sum_{g\in G}[\hh_i^g, \hh_i^g, x_i])= \sum_{i=1}^s(\sum_{g\in G}\a(\hh_i^g, \hh_i^g)x_i=0.
$$
Then $x_i\in \hh,$ which is  a contradiction. Thus
$$
\LL_i=\hh_i\oplus\bigoplus_{g\in G}(\bigoplus_{\a\in\Pi_i^g}\LL_{\a, i}^g)=\hh_i\oplus\bigoplus_{\a\in\Pi_i}\LL_{\a, i},
$$
where
\begin{equation}\label{43.}
\Pi_i^g=\{\a|_{\hh_i\times\hh_i}~:~ \LL_{\a, i}^g=\LL_i\cap\LL_\a^g\neq0\},~~\Pi_i=\bigcup_{g\in G}\Pi_i^g.
\end{equation}
Note that $\Pi=\bigcup_{i=1}^s\Pi_i.$
By Definition \ref{ideal}, we have $\rho(\LL_i, \LL)(A \LL)\subset \LL_i$ and $\rho(\LL_i, \LL_i)(A \LL_i)\subset \LL_i,~~i=1, 2, 3, ..., s.$ Since  $Z_{\LL}(A) = 0,$ we also have $\rho(\LL_i, \LL_j)A=0$ when $i\neq J.$ Now, we  can be written as
$$
A_i=\bigoplus_{g\in G}(A_{0, i}^g\oplus(\bigoplus_{\lam\in\Lam_i^g}A_{\lam, i}^g)),
$$
where for any $g\in G$, 
\begin{eqnarray*}
A_{0, i}^g&=&\{a\in A^g~:~\rho(\hh_i^g, \hh_i^g)a=0\}\\
A_{\lam, i}^g&=&\{a\in A^g~:~\rho(h,h')a=\lam(h, h')a,~~\forall h, h'\in\hh_i^g,~~\lam\in\Lam,~~\lam(\hh_i^g, \hh_i^g)\neq 0\},
\end{eqnarray*}
and
\begin{equation}\label{44.}
\Lam_i^g=\{\lam|_{\hh_i^g\times\hh_i^g}~:~\lam\in\Lam,~~\lam(\hh_i^g, \hh_i^g)\neq 0\},~~\Lam_i=\bigcup_{g\in G}\Lam_i^g.
\end{equation}
Note that $\Lam=\bigcup_{i=1}^s\Lam_i.$
Therefore, $(\LL_i, A, \rho|_{\LL_i\times\LL_i})$ is a split $3-$Lie-Rinehart color algebra with  a splitting Cartan subalgebra $\hh_i$ and the root system $\Pi_i$ defined by (\ref{43.}), and the weight system $\Lam_i$ defined by (\ref{44.}) associated to $\hh_i.$

Conversely, suppose that $(\LL_i, A, \rho|_{\LL_i\times\LL_i}), 1 \leq i \leq s$ are split $3-$Lie-Rinehart color algebras, and $\LL_i=\hh_i\oplus\bigoplus_{g\in G}(\bigoplus_{\a\in\Pi_i^g}\LL_{\a, i}^g)=\hh_i\oplus\bigoplus_{\a\in\Pi_i}\LL_{\a, i},$ with a splitting Cartan subalgebra $\hh_i,$ the root system $\Pi_i,$ and the weight system $\Lam_i$ associated to $\hh_i,$ respectively. By taking $\hh =\bigoplus_{i=1}^s\hh_i$ being an abelian subalgebra of $\LL=\bigoplus_{i=1}^s\LL_i.$ Now, for any $\a\in\Pi_i$ with  $\a|_{\hh_i\times\hh_i}\neq 0,$ we extend $\a$ to $\a~:~\hh\times\hh\longrightarrow\bbbk$ by
 \begin{equation*}
	\a(h, h')~:=\left\{
	\begin{array}{ll}
		0~~~~~~~&;~~h\notin\hh_i ~\hbox{~or~}~h'\notin\hh_i\\
		\a(h, h')~~~~~~~~~~&;~~~~h, h'\in\hh_i.
	\end{array}
	\right.
\end{equation*}
We get that $\LL_\a=\LL_{\a, i}\neq 0$ and $\a\in (\hh\times\hh)^*\neq 0.$ Therefore, $\hh =\bigoplus_{i=1}^s\hh_i$ is a splitting Cartan subalgebra of $\LL$ with root system $\Pi=\bigcup_{i=1}^s\Pi_i$ and $\LL=\hh\oplus\bigoplus_{\a\in\Pi}\LL_\a.$ Thanks to  Definition \ref{ideal} and  $Z_{\LL}(A) = 0$, $\rho(\LL_i, \LL_j)A=0$ for $i\neq J,$ and $\rho(\LL_i, \LL_i)(A \LL_i)\subset \LL_i,~~1 \leq i \leq s.$  Therefore, $\rho(\hh_i
, \hh_i)A\LL_i = 0, 1 \leq i, j \leq s.$
By a complete similar discussion to the above, for any $\lam\in\Lam_i,$
 extending $\lam$ to $\lam~:~\hh\times\hh\longrightarrow\bbbk.$ we get that $A_\lam = A_i\neq 0,
$ and $\lam\in(\hh\times\hh)^*\neq 0.$ Therefore, $\Lam=\bigcup_{i=1}^s\Lam_i\subset (\hh\times\hh)^*\neq 0 $ and $A=A_0\oplus\bigoplus_{\lam\in\Lam}A_\lam,$ where $A_0=\bigcap_{i=1}^s A_{0, i}.$ We get the result.\qed

\begin{DEF}  A split $3-$Lie-Rinehart color algebra $(\LL, A)$ is tight if $Z_\rho(\LL) = Ann(A) = 0,
	AA = A, A\LL = \LL $ and
\begin{eqnarray*}
\hh=(\sum_{\substack{\b\in\Pi, -\b\in\Lam\\ g, g'\in G}} A_{-\b}^g\LL_\b^{g'})
+(\sum_{\substack{\b, \gamma, \mu\in\Pi\\ \b+\gamma+\mu=0\\ g, g', g''\in G}}[\LL_\b^g,
\LL_{\gamma}^{g'}, \LL_\mu^{g''}]),\\
A_0=(\sum_{\substack{\mu\in\Lam\\  g\in G}} A_\mu^{g} A_{-\mu}^{-g})
+(\sum_{\substack{\a, \b\in\Pi, \mu\in\Lam\\ \a+\b+\mu=0\\ g, g', g''\in G}}\rho(\LL_\a^g,
\LL_\b^{g'}) A_\mu^{g''}).
\end{eqnarray*}
\end{DEF} 

\begin{rem} If $(\LL, A)$ is a tight  split $3-$Lie-Rinehart color algebra then it follows from Proposition \ref{2.17} and Corolary \ref{4.19} that 
$$
\LL=\bigoplus_{[\a]\in\Pi/\sim}I_{[\a]},~~ A=\bigoplus_{[\lam]\in\Lam/\approx}\aa_{[\lam]},
$$	
with any $I_{[\a]}$ an ideal of $\LL$ satisfying $$
[I_{[\a]}, I_{[\b]}, I_{[\gamma]}]=0,~~~\hbox{and~~}~~~~[I_{[\a]}, I_{[\a]}, I_{[\b]}]=0,
$$
whenever $[\a], [\b], [\gamma]\in\Pi/\sim$ be different from each other, and any $\aa_{[\lam]}$ an ideal of $A$  satisfying  $\aa_{[\lam]} \aa_{[\mu]}=0,$ when $[\lam]\neq [\mu].$
\end{rem}

\begin{pro}\label{5.10}   Let $(\LL, A)$ be a tight split $3-$Lie-Rinehart color algebra. Then for any $[\a]\in\Pi/\sim$
there exists a unique $[\lam]\in\Lam/\approx$ such that $ \aa_{[\lam]} I_{[\a]} \neq 0.$ 
\end{pro}
\noindent {\bf Proof.} It can be analogous Proposition 4.2 in
\cite{ABCS} .\qed

\begin{thm}\label{5.19}  Let $(\LL, A)$ be a tight split $3-$Lie-Rinehart color algebra. Then
$$
\LL=\bigoplus_{g\in G}\bigoplus_{i\in I}\LL_i^g=\bigoplus_{i\in I}\LL_i,~~A=\bigoplus_{h\in G}\bigoplus_{j\in J}A_j^h=\bigoplus_{j\in J} A_j,
$$
where any $\LL_i=\bigoplus_{g\in G}\LL_i^g$ is a non-zero graded ideal of $\LL$ satisfying $[\LL_{i_1}^{g_1}, \LL_{i_2}^{g_2}, \LL_{i_3}^{g_3}]=0,$ when $i_1, i_2, i_3\in I, g_1, g_2, g_3\in G$ be different from each other, and any $ A_j=\bigoplus_{h\in G}A_j^h$ is a graded ideal of $A$ such that $A_{j_1}^{h_1} A_{j_2}^{h_2}=0$ when $(j_1, h_1)\neq(j_2, h_2).$ Moreover, both decompositions satisfy that for any $(i, g)\in I\times G$ there exists a unique $(j, h)\in J\times G$ such that 
$$
A_j^h \LL_i^g\neq 0. 
$$
Fortheremore, any $(\LL_i=\bigoplus_{g\in G}\LL_i^g, A_j=\bigoplus_{h\in G}A_j^h, \rho|_{\LL_i\times\LL_i})$ is a split $3-$Lie-Rinehart color algebra. 
\end{thm}
\noindent {\bf Proof.}  Proposition \ref{5.10} shows that $\LL_i$
is an $A_j-$module. Hence we can get the rsults of theorem. For the final result of theorem see Proposition \ref{5.4}.\qed

\section{the  simple components of  split $3$-Lie-Rinehart  color algebras} 
\setcounter{equation}{0}\

In this section we focus on the simplicity of  split $3$-Lie-Rinehart  color algebra
$(\LL, A)$ by centering our attention in those of maximal length. From now on, we will suppose $\Pi$ and $\Lam$ are symmetric.

Let us introduce the concepts of root-multiplicativity and
maximal length in the framework  of  split $3$-Lie-Rinehart  color algebra, in a similar way to the ones for  split
Lie-Rinehart algebra in \cite{ABCS}.

\begin{DEF}\label{mal} A split  $3$-Lie-Rinehart  color algebra
$(\LL, A)$  is called  {\em root-multiplicative} if for any $\a, \b, \gamma\in\Pi$ and $\lam, \mu\in\Lam$ the following conditions hold 
\begin{itemize}
	\item[-] If  $\a+\b+\gamma\in\Pi^g, g\in G$\ then $[\LL_\a^g, \LL_\b^{g'}, \LL_\gamma^{g''}]\neq 0.$
	\item[-]  If  $\lam+\a\in\Pi^g, g\in G$ then $A_\lam^g\LL_\a^{g'}\neq 0.$
	\item[-]  If  $\lam+\mu\in\Lam^g, g\in G$ then $A_\lam^g A_\mu^{g'}\neq 0.$
\end{itemize}
\end{DEF}

\begin{DEF} A split  $3$-Lie-Rinehart  color algebra
	$(\LL, A)$  is called of {\em maximal length} if  for
any $\a\in\Pi^g, g\in G$ and  $\lam\in\Lam^h, h\in G$  we have $\dim
\LL^{g}_{\a}=\dim A_\lam^h=1.$ 
\end{DEF}

\begin{rem} If $(\LL, A)$ is a split  $3$-Lie-Rinehart  color algebra such that $\LL$ and $A$ are simple algebras then $Z_\rho(\LL) = Ann(A) = \{0\}$. Also as consequence of Theorem \ref{main1}-(2) and Theorem \ref{main1'}-(2) we get that all of the non-zero roots
in $\Pi$ are connected, that all of the non-zero weights in $\Lam$ are also connected and that
\begin{eqnarray*}
\hh=(\sum_{\substack{\b\in\Pi\cap\Lam\\ g, g'\in G}} A_{-\b}^g\LL_\b^{g'})
+(\sum_{\substack{\b, \gamma, \mu\in\Pi\\ \b+\gamma+\mu=0\\ g, g', g''\in G}}[\LL_\b^g,
\LL_{\gamma}^{g'}, \LL_\mu^{g''}]),\\
A_0=(\sum_{\substack{\mu\in\Lam\\  g\in G}} A_\mu^{g} A_{-\mu}^{-g})
+(\sum_{\substack{\a, \b\in\Pi, \mu\in\Lam\\ \a+\b+\mu=0\\ g, g', g''\in G}}\rho(\LL_\a^g,
\LL_\b^{g'}) A_\mu^{g''}).
\end{eqnarray*}
From here, the conditions for $(\LL, A)$ of being tight 
together with the ones of having $\Pi$ and $\Lam$ all of their elements connected, are necessary
conditions to get a characterization of the simplicity of the algebras $\LL$ and $A.$ Actually, we are going to shwo that 
under the hypothesis of being $(\LL, A)$ of maximal length and root-multiplicative, these are
also sufficient conditions. 
\end{rem}

\begin{pro}\label{6.4} Let $(\LL, A)$ be a tight  split $3$-Lie-Rinehart  color algebra of maximal length and root-multiplicative. If $\LL$ has
	all of its nonzero roots connected,  then either $\LL$ is simple or $\LL=I\oplus I'$ where  $I$ and $I'$ are simple ideals of $\LL.$ 
\end{pro}	
\noindent {\bf Proof.} Consider any nonzero ideal $I$ of $\LL,$
	by Lemma \ref{5.3}, $I\nsubseteq\hh^g, \forall g\in G$ and by Eq. (\ref{301}) we can write
\begin{equation*}
		I=\bigoplus_{g\in G}((I\cap\hh^g)\oplus(\bigoplus_{\a\in\Pi^g_I}(I\cap\LL_\a^g))),
\end{equation*}
where $\Pi^g_I :=\{\a\in\Pi : I_\a^g=I\cap\LL_\a^g\neq0\}\subset\Pi^g$ for all $g\in G$  and some
$\Pi^g_I\neq \phi.$ Then we can write $	I=\bigoplus_{g\in G}((I\cap\hh^g)\oplus(\bigoplus_{\a\in\Pi^g_I}I_\a^g)).$  Let us
distinguish two cases.
	
{\bf Case 1.} Suppose there exists $\a_0\in\Pi^g_I$ such that $-\a_0\in\Pi^g_I$  for all $g\in G.$ Then $0\neq I_{\a_0}^g$ and  by the maximal length of $(\LL, A)$  that
\begin{equation}\label{302.}
0\neq\LL_{\a_0}^g\subset I.
\end{equation}

Now, let us take some $\b\in\Pi$ satisfying $\b\notin\{\a_0, -\a_0\}.$ Since $\a_0$ and $\b$
	are connected, we have a connection  $\{\b_1, \b_2,...,
	\b_{2k+1}\},~~k\geq 2$ from $\a_0$ to
	$\b$ satisfying the following conditions;\\

\begin{itemize}
	\item[(1)] $\b_1=\pm\a_0.$
	\item[(2)] $\b_1+(\b_2+\b_3)\in\pm\pi,\\
	\b_1+ (\b_2+( \b_3+(\b_4+\b_5)))\in\pm\Pi,\\
	...\\
	\b_i+ \sum_{j=1}^i(\b_{2j}+ \b_{2j+1})\in\pm\Pi,~~0< i< k.$
	
	\item[(3)] $\b_1+ \sum_{j=1}^k(\b_{2j}+ \b_{2j+1})\in\{\b, -\b\}.$
\end{itemize}	
Consider $\b_1, \b_2, \b_3\in\Pi\cup\Lam$ and $\b_2+\b_3\in\Pi,$ since
	$\b_1\in\Pi$ there exist $g_1\in\Lam$ such that
	$\LL_{\b_1}^{g_1}\neq 0.$ From here, the root-multiplicativity
	and maximal length of $\LL$ allow us to get
	$$
	0\neq[\LL_{\b_1}^{g_1},
	\LL_{\b_2}^{g_2}, \LL_{\b_3}^{g_3}]=\LL_{\b_1+\b_2+\b_3}^{g_1+g_2+g_3}.
	$$
	Since $0\neq\LL_{\b_1}^{g_1}\subset I$ as consequence of Eq.
	(\ref{302}) we have
	$$
	0\neq\LL_{\b_1+\b_2+\b_3}^{g_1+g_2+g_3}\subset I.
	$$
	We can argue in a similar way from $\b_1+(\b_2+
	\b_3),~~\b_4$ and $\b_1+((\b_2
	+\b_3)+\b_4) \in\Pi$ to get
	$$
	0\neq\LL_{\b_1+\b_2+\b_3+\b_4}^{g_1+g_2+g_3+g_4}\subset I, ~~ \hbox{~for~~ some}~~
	g_4\in G.
	$$
	We can follow this process with the connection $\{\b_1, \b_2,...,
	\b_{2k+1}\}$ we obtain that
	$$
	0\neq\LL_{\b_1+ \b_2+\b_3+
	...+\b_{2k+1}}^{h}\subset I, ~~
	\hbox{~for~~ some}~~ h\in G.
	$$
	Thus we have shown that 
\begin{equation}\label{305.}
\hbox{for ~~any~~} \b\in\Pi,~ \hbox{we~~ have~~ that~}~ 0\neq \LL_{\varepsilon\b}^h\subset I\hbox{~~for~~some~~}~~\varepsilon\in\{\pm 1\}.
\end{equation}	
Since $-\b\in\Pi_I^g$ we have $\{-\b_1, -\b_2, ..., -\b_{2k+1}\}$ is a connection from $-\b$ to $\a_0$ satisfying
$$
\{-\b_1, -\b_2, ..., -\b_{2k+1}\}\in\{\b, -\b\}.
$$
By arguing as above we get,
\begin{equation}\label{306.}
0\neq \LL_{-\varepsilon\b}^g\subset I,	
\end{equation}	
and so $\bigcup_{g\in G}\Pi_I^g=\Pi.$ From the fact 
$$
\hh=(\sum_{\substack{\b\in\Pi\cap\Lam\\ g, g'\in G}} A_{-\b}^g\LL_\b^{g'})
+(\sum_{\substack{\b, \gamma, \mu\in\Pi\\ \b+\gamma+\mu=0\\ g, g', g''\in G}}[\LL_\b^g,
\LL_{\gamma}^{g'}, \LL_\mu^{g''}]),
$$
imply that
	\begin{equation}\label{330.}
		\hh\subset I.
	\end{equation}
From Eqs. (\ref{302.})-(\ref{330.}), ) we obtain $\LL\subset I,$  and so $\LL$ is simple.

{\bf Case 2.} In the second case, suppose that for any $\a_0\in\Pi^g_I$ we have that $-\a_0\notin\Pi^g_I,$ for all $g\in G.$ Observe that by
arguing as in the previous case we can write
\begin{equation}\label{340.}
\Pi^g=\Pi_I^g\cup-\Pi_I^g,~~~\forall g\in G,
\end{equation}
where $-\Pi_I^g=\{-\a~:~ \a\in\Pi^g_I\}.$ Denote by
$$
I'~:=\bigoplus_{g\in G}(\sum_{\substack{\a\in\Lam^g\\ -\a\in -\Pi_I^g}}A_\a^g\LL_{-\a}^{-g}\oplus(\bigoplus_{-\a\in-\Pi_I^g}\LL_{-\a}^{-g})).
$$
We are giong to show that $I'$ is a $3-$Lie color ideal of $\LL.$ First, we will show that $I'$ is a Lie ideal of $\LL.$ Taking into account Eq. (\ref{00}) and the fact that $\hh$ is abelian we have
\begin{eqnarray}\label{350}
\nonumber[\LL, \LL, I']&=&[\hh\oplus\bigoplus_{\substack{\b\in\Pi^h\\ h\in G}}\LL^h_\b, \hh\oplus\bigoplus_{\substack{\gamma\in\Pi^l\\ l\in G}}\LL^l_\gamma, \sum_{\substack{\a\in\Lam^g\\ -\a\in -\Pi_I^g\\ g\in G}}A_\a^g\LL_{-\a}^{-g}\oplus(\bigoplus_{\substack{-\a\in-\Pi_I^g\\ g\in G}}\LL_{-\a}^{-g}))]\\
&\subset&\bigoplus_{\substack{-\a\in-\Pi_I^g\\ g\in G}}\LL_{-\a}^{-g}+[\bigoplus_{\substack{\b\in\Pi^h\\ h\in G}}\LL^h_\b, \hh,  \sum_{\substack{\a\in\Lam^g\\ -\a\in -\Pi_I^g\\ g\in G}}A_\a^g\LL_{-\a}^{-g} ]\\
\nonumber&+&[\bigoplus_{\substack{\b\in\Pi^h\\ h\in G}}\LL^h_\b, \hh, \bigoplus_{\substack{-\a\in-\Pi_I^g\\ g\in G}}\LL_{-\a}^{-g}]+[\bigoplus_{\substack{\b\in\Pi^h\\ h\in G}}\LL^h_\b, \bigoplus_{\substack{\gamma\in\Pi^l\\ l\in G}}\LL^l_\gamma, \sum_{\substack{\a\in\Lam^g\\ -\a\in -\Pi_I^g\\ g\in G}}A_\a^g\LL_{-\a}^{-g}]\\
\nonumber&+&[\bigoplus_{\substack{\b\in\Pi^h\\ h\in G}}\LL^g_\a, \bigoplus_{\substack{\gamma\in\Pi^l\\ l\in G}}\LL^l_\gamma, \bigoplus_{\substack{-\a\in-\Lam^g\\ g\in G}}\LL^{-g}_{-\a}].
\end{eqnarray}

Consider the second summand in Eq. (\ref{350}). If some $[\LL_\b^h, \hh^l, A_\a^g\LL_{-\a}^{-g}]\neq 0$ and  $\b=-\a,$ we have that $[\LL_{-\a}^h, \hh^l, A_\a^g\LL_{-\a}^{-g}]\subset\LL_{-\a}^{h+l}\subset I'.$ In case $\b=\a,$ since $I$ is an ideal of $\LL,$ we have $-\a\notin \Pi_I^g$ implies $[\LL_{-\a}^h, \hh^l, A_{-\a}^g\LL_{-\a}^{-g}]=0,$ by symmetry of $\Pi$ and maximal lenght of $\LL,$ we have $[\LL_\a^h, \hh^l, A_\a^g\LL_{-\a}^{-g}]=0.$ Suppose that $\b\notin\{\a, -\a\}.$ As $[\LL_\b^h, \hh^l, A_\a^g\LL_{-\a}^{-g}]\neq 0,$ By Eq. (\ref{4}), either $ A_\a^g[\LL_\b^h, \hh^l,\LL_{-\a}^{-g}]\neq 0$ or $\rho(\LL_\b^h,\hh^l)(A_\a^g)\LL_{-\a}^{-g}\neq 0.$  By the maximal length of $\LL,$  either $ A_\a^g[\LL_\b^h, \hh^l,\LL_{-\a}^{-g}]=\LL_\b^{h+l}$ or $\rho(\LL_\b^h,\hh^l)(A_\a^g)\LL_{-\a}^{-g}=\LL_\b^{h+l}.$ In both cases, since $\a\in\Pi_I^g,$ by root-multiplicativity, we have $\LL_{-\b}\subset I$ and therefore $-\b\in\Pi_I^g,$ that is $\LL_\b^{h+l}\subset I'.$ Thus  
\begin{equation}\label{360}
[\bigoplus_{\substack{\b\in\Pi^h\\ h\in G}}\LL^h_\b, \hh,  \sum_{\substack{\a\in\Lam^g\\ -\a\in -\Pi_I^g\\ g\in G}}A_\a^g\LL_{-\a}^{-g}]\subset I'. 
\end{equation}
A similar discussion as above, one can show that
\begin{equation}\label{361}
[\bigoplus_{\substack{\b\in\Pi^h\\ h\in G}}\LL^h_\b, \bigoplus_{\substack{\gamma\in\Pi^l\\ l\in G}}\LL^l_\gamma, \sum_{\substack{\a\in\Lam^g\\ -\a\in -\Pi_I^g\\ g\in G}}A_\a^g\LL_{-\a}^{-g}]\subset I'. 
\end{equation}

Now, if we consider the third summand in Eq. (\ref{350}) and some $[\LL_\b^h, \hh^l, \LL_{-\a}^{-g}]\neq 0.$ Then we have $[\LL_\b^h, \hh^l, \LL_{-\a}^{-g}]=\LL_{-\a+\b}^{-g+h+l}.$ If $\b\neq\a,$ thanks to $\a\in\Pi_I^g$ and the root-multiplicativity give us $[\LL_\a^g, \hh^{-l} \LL_{-\b}^{-h}]=\LL_{\a-\b}^{g-h-l}\subset I.$ Hence $-\a+\b\in-\Pi_I^{-g+h+l}$ and then $\LL_{-\a+\b}^{-g+h+l}\subset I'.$ Now, if $\b=\a,$ in case $[\LL_\a^h, \hh^l, \LL_{-\a}^{-g}]\neq 0,$ we get  $[\LL_\a^h, \hh^l, \LL_{-\a}^{-g}]\subset I,$ because of $\a\in\Pi_I^g.$ Thus 
$$
\LL_{-\a}^{-g}=[[\LL_\a^h, \hh^l, \LL_{-\a}^{-g}], \hh^{-l}, \LL_{-\a}^{-h}]\subset I.
$$
From here, we deduce $-\a, \a\in\Pi_I^g,$ a contradiction with Eq. (\ref{340.}). Thus
\begin{equation}\label{362}
[\bigoplus_{\substack{\b\in\Pi^h\\ h\in G}}\LL^h_\b, \hh, \bigoplus_{\substack{-\a\in-\Pi_I^g\\ g\in G}}\LL_{-\a}^{-g}]\subset I'.
\end{equation}
Again a similar discussion as above, one can show that
\begin{equation}\label{363}
[\bigoplus_{\substack{\b\in\Pi^h\\ h\in G}}\LL^g_\a, \bigoplus_{\substack{\gamma\in\Pi^l\\ l\in G}}\LL^l_\gamma, \bigoplus_{\substack{-\a\in-\Lam^g\\ g\in G}}\LL^{-g}_{-\a}]\subset I'.
\end{equation}
From Eqs. (\ref{360})-(\ref{363}), we conclud that $I'$ is a $3-$Lie color ideal of $\LL.$
Second, we will check $A I'\subset I'.$  Taking into account Eq. (\ref{00.}),  we have

\begin{eqnarray}\label{370}
\nonumber A I'&=&(A^0\oplus(\bigoplus_{h\in G}\bigoplus_{\lam\in\Lam^h}A^h_\lam))(\bigoplus_{g\in G}(\sum_{\substack{\a\in\Lam^g\\ -\a\in -\Pi_I^g}}A_\a^g\LL_{-\a}^{-g}\oplus(\bigoplus_{-\a\in-\Pi_I^g}\LL_{-\a}^{-g})))\\
&\subset&I'+(\bigoplus_{\substack{\lam\in\Lam^h\\ h\in G}}A^h_\lam)(\sum_{\substack{\a\in\Lam^g\\ -\a\in -\Pi_I^g\\ g\in G}}A_\a^g\LL_{-\a}^{-g})+(\bigoplus_{\substack{\lam\in\Lam^h\\ h\in G}}A^h_\lam)(\bigoplus_{\substack{-\a\in-\Pi_I^g\\ g\in G}}\LL_{-\a}^{-g})
\end{eqnarray}
Consider the third summand in (\ref{370}) and suppose that $A_\lam^h \LL_{-\a}^{-g}\neq 0$ for some $\lam\in\Lam^h,~-\a\in\Pi_I^g.$ If $\a-\lam\in\Pi_I^{g-h},$ then by the root-multiplicativity of $\LL$ we get $A_{-\lam}^{-h}\LL_\a^g\neq 0.$  Now by the maximal length of $\LL$ and the fact $\a\in\Pi_I^g,$ we conclud that $A_{-\lam}^{-h}\LL_\a^g=\LL_{\a-\lam}^{g-h}\subset I.$ Therefore $-\a+\lam\in\Pi_I^{-g+h}$ which is a contradiction. Hence $\a-\lam\in-\Pi_I^{g-h},$ and so $A_\lam^h \LL_{-\a}^{-g}\subset I'.$ Therefore,
\begin{equation}\label{363}
(\bigoplus_{\substack{\lam\in\Lam^h\\ h\in G}}A^h_\lam)(\bigoplus_{\substack{-\a\in-\Pi_I^g\\ g\in G}}\LL_{-\a}^{-g})\subset I'.
\end{equation}
We can argue as above with the second summand in (\ref{370}) so as to conclude that
\begin{equation}\label{375.}
(\bigoplus_{\substack{\lam\in\Lam^h\\ h\in G}}A^h_\lam)(\sum_{\substack{\a\in\Lam^g\\ -\a\in -\Pi_I^g\\ g\in G}}A_\a^g\LL_{-\a}^{-g})\subset I'.
\end{equation}
From Eqs. (\ref{363}) and (\ref{375.}) we get $A I'\subset I'.$
Finally, let us check $\rho(I', I')(A)\LL\subset I'.$ In fact by Eq. (\ref{4}) we have
$$
\rho(I', I')(A)\LL\subset [I', I', A\LL]+A[I', I', \LL]
$$
Tanks to $I'$ is a $3-$Lie color ideal we get the result. 

Summarizing a discussion of above, we conclude that $I'$ is an ideal of the split $3-$Lie-Rinehart color algebra $(\LL, A).$

Next, by Eq. (\ref{340.}) we get $\sum_{\substack{\b, \gamma, \mu\in\Pi\\ \b+\gamma+\mu=0\\ g, g', g''\in G}}[\LL_\b^g,
\LL_{\gamma}^{g'}, \LL_\mu^{g''}]=0,$ so by hypothesis must have 
$$
\hh=\sum_{\substack{\a\in\Pi_I^g\\ -\a\in\Lam_I^g\\ g\in G}} A_{-\a}^{-g}\LL_\a^g\oplus\sum_{\substack{-\a\in-\Pi_I^g\\ \a\in\Lam_I^g\\ g\in G}} A_\a^g\LL_{-\a}^{-g}.
$$
For direct character, take
$$
0\neq h, h'\in\sum_{\substack{\a\in\Pi_I^g\\ -\a\in\Lam_I^g\\ g\in G}} A_{-\a}^{-g}\LL_\a^g\cap\sum_{\substack{-\a\in-\Pi_I^g\\ \a\in\Lam_I^g\\ g\in G}} A_\a^g\LL_{-\a}^{-g}.
$$
Taking into account $Z_\rho(\LL) = \{0\}$ and $\LL$ is split, there exists $0\neq x\in\LL_\b^h, \b\in\Pi^h,$ such
that $[h, h', x]\neq 0,$ being then $x\in I\cap I'=\{0\},$ a contradiction. Hence  the sum is direct. Taking into account the above observation and Eq. (\ref{340.}) we
have
$$
\LL=I\oplus I'.
$$
Finally, we can proceed with $I$ and $I'$ as we did for $\LL$ in the first case of the proof to
conclude that $I$ and $I'$ are simple ideals, which completes the proof of the theorem.\qed

In a similar way to Proposition \ref{6.4} one can prove the next result;

\begin{pro}\label{6.4.} Let $(\LL, A)$ be a tight  split $3$-Lie-Rinehart  color algebra of maximal length and root-multiplicative. If $A$ has
all of its nonzero ero weights connected,  then either $A$ is simple or $A=J\oplus J'$ where  $J$ and $J'$ are simple ideals of $A.$\qed 
\end{pro}	

Now, we are ready to state our main result;

\begin{thm}\label{final}  Let $(\LL, A)$ be a tight split $3-$Lie-Rinehart color algebra of maximal length, root-multiplicative, with symmetric roots and weight systems in such a way that the root system $\Pi$ has all its elements connected and the weight system $\Lam$ has all its elements connected.Then
$$
\LL=\bigoplus_{g\in G}\bigoplus_{i\in I}\LL_i^g=\bigoplus_{i\in I}\LL_i,~~A=\bigoplus_{h\in G}\bigoplus_{j\in J}A_j^h=\bigoplus_{j\in J} A_j,
$$
where any $\LL_i=\bigoplus_{g\in G}\LL_i^g$ is a simple  graded ideal of $\LL$ satisfying $[\LL_{i_1}^{g_1}, \LL_{i_2}^{g_2}, \LL_{i_3}^{g_3}]=0,$ when $i_1, i_2, i_3\in I, g_1, g_2, g_3\in G$ be different from each other, and any $ A_j=\bigoplus_{h\in G}A_j^h$ is a simple graded ideal of $A$ such that $A_{j_1}^{h_1} A_{j_2}^{h_2}=0$ when $(j_1, h_1)\neq(j_2, h_2).$ Moreover, both decompositions satisfy that for any $(i, g)\in I\times G$ there exists a unique $(j, h)\in J\times G$ such that 
$$
A_j^h \LL_i^g\neq 0. 
$$
Fortheremore, any $(\LL_i=\bigoplus_{g\in G}\LL_i^g, A_j=\bigoplus_{h\in G}A_j^h, \rho|_{\LL_i\times\LL_i})$ is a split $3-$Lie-Rinehart color algebra. 
\end{thm}
\noindent {\bf Proof.} By Theorem \ref{5.19}  we can can write
$$
\LL=\bigoplus_{[\a]\in\Pi/\sim}I_{[\a]},
$$
with any $I_{[\a]}$ an ideal of $\LL,$ being each $I_{[\a]}$ a split $3-$Lie-Rinehart color algebra having as root system $[\a].$  Also we can write $A$ as the direct sum of the ideals
$$
 A=\bigoplus_{[\lam]\in\Lam/\approx}\aa_{[\lam]},
$$
in such a way that any $\aa_{[\lam]}$ has as weight system $[\lam],$ for any $[\a]\in\Pi/\sim$
there exists a unique $[\lam]\in\Lam/\approx$ such that $ \aa_{[\lam]} I_{[\a]} \neq 0$ and being $( I_{[\a]}, \aa_{[\lam]})$ a split $3-$Lie-Rinehart color algebra.

In order to apply Proposition \ref{6.4} and Proposition \ref{6.4.} to each $( I_{[\a]}, \aa_{[\lam]}),$  we previously
have to observe that the root-multiplicativity of  $( I_{[\a]}, \aa_{[\lam]}),$ Proposition \ref{12345}  and Theorem \ref{main2}
show that $[\a]$ and $[\lambda]$ have, respectively, all of their elements $[\a,~[\lambda]-]$connected. That is,
connected through connections contained in $[\a]$  and $[\lambda].$  Any of the $( I_{[\a]}, \aa_{[\lam]})$ is root-multiplicative as consequence of the root-multiplicativity of $(\LL, A).$  Clearly  $( I_{[\a]}, \aa_{[\lam]})$
is of maximal length and tight, last fact consequence of tightness of $(\LL, A),$ Proposition \ref{6.4} and Proposition \ref{6.4.}. So we can apply Proposition \ref{6.4} and Proposition \ref{6.4.} to each $( I_{[\a]}, \aa_{[\lam]})$
so as to conclude that any $ I_{[\a]}$
is either simple or the direct sum of simple
ideals $ I_{[\a]}=J\oplus J',$ and that any $\aa_{[\lam]}$
is either simple or the direct sum of simple ideals $\aa_{[\lam]}=B\oplus B'.$ From here, it is clear that by writing $I_i=J\oplus J'$ and $\aa_j=B\oplus B'$ if $I_i$ or $\aa_j$ are not, respectively, simple, then Theorem \ref{5.19} allows as to assert that the resulting
decomposition satisfies the assertions of the theorem. \qed


\end{document}